\documentclass{article}

\usepackage[a4paper,top=3cm,bottom=3cm,left=2cm,right=2cm]{geometry}
\usepackage{amsmath}
\usepackage{amssymb}
\usepackage{amsthm}
\usepackage{color}
\usepackage{multicol}
\usepackage{graphicx}
\usepackage{algpseudocode} 
\usepackage{algorithm} 
\usepackage{picture}
\usepackage{tikz}
\usepackage{accents}
\usepackage{cases}
\usepackage{float}
\usepackage{graphicx}
\usepackage{epstopdf}
\usepackage{subcaption} 
\usepackage[figuresright]{rotating}

\renewcommand{\d}{\mathrm{d}}
\newcommand{\vect}[1]{\boldsymbol{#1}}

\renewcommand{\d}{\mathrm{d}}
\newcommand{\dt}{\mathrm{dt}}

\renewcommand{\S}{\mathcal{S}}

\newcommand{\vecx}{\boldsymbol{x}}

\begin{document}

\title{High-Order Spline Upwind for Space-Time Isogeometric Analysis}
\author{G. Loli$^1$, G. Sangalli$^{1,2}$ and P. Tesini$^{3,1}$}
\date{}
\maketitle

\begin{center}
\small{$^1$ \ Universit\`a di Pavia, Dipartimento di Matematica ``F. Casorati'' \\ 
Via A. Ferrata 5, 27100 Pavia, Italy. \\ 
\{gabriele.loli, giancarlo.sangalli\}@unipv.it} 
\end{center}
\begin{center}
\small{$^2$ \ IMATI-CNR ``Enrico Magenes'',  Pavia, Italy.}
\end{center}
\begin{center}
\small{$^3$ \ Universit\`a degli Studi di Milano-Bicocca \\
Piazza dell'Ateneo Nuovo 1, 20126 Milano, Italy. \\ 
p.tesini@campus.unimib.it}
\end{center} 

\begin{abstract}
 We propose an innovative isogeometric space-time method for the heat
  equation, with  smooth splines approximation in both space
  and time. To enhance the stability of the method we add a stabilizing term, based on a linear combination of
  high-order artificial diffusions. This term  is designed in order to make the linear
  system lower block-triangular, that is, lower triangular with
  respect to time. In order to keep optimal
  accuracy, the stabilization terms are further weighted  in terms
   of the residual. Through a series of numerical experiments, we validate
  the method's capability, showcasing its stability and accuracy.
\end{abstract}

\begin{flushleft}
\end{flushleft}

\section{Introduction}

Isogeometric Analysis (IgA), introduced in \cite{Hughes2005} (see also
the book \cite{Cottrell2009}), is an evolution of the classical
finite element method. IgA  uses spline  functions, or their
generalizations, both to  represent the computational domain and to
approximate the solution of the partial differential equation that
models the problem of interest, in order to facilitate the interoperability between
computer aided design (CAD) and numerical simulation. At the same
time,  IgA benefits from the properties of smooth splines, such as higher accuracy when compared to $C^0$ piecewise polynomials (see e.g. \cite{Evans2009, Bressan2018}).

The idea of using finite elements in the  space-time domain comes from
\cite{Fried1969, Oden1969, Argyris1969} and was then developed for various problems as heat transfer
\cite{Bruch1974}, advection-diffusion \cite{Nguyen1984} and
elastodynamics \cite{Hughes1988}.  The  mathematical theory  of
space-time Galerkin methods has been developed in recent works, for example
\cite{Schwab2009, Steinbach2015} .

Space-time formulations in  IgA
provides an additional opportunity, that is, to exploit the properties of smooth splines in
time as well, as proposed in 
\cite{Takizawa2014, Langer2016}. In particular, \cite{Langer2016}
develops a stabilized IgA  of the heat equation.  In
\cite{Montardini2018, Loli2020} the authors have proposed
preconditioners and solvers, while in \cite{Saade2021} a continuous
space-time IgA formulation has been applied to linear and non-linear
elastodynamics. The use of smooth splines with respect to time  poses interesting
challenges as well.

A challenge with space-time
formulations concerns the
 causality principle. While the sequentiality of  discontinuous Galerkin in time
guarantees causality,  this is not the case for
 Galerkin with  smooth
 spline approximation in time. The lack of
 causality generates further unphysical behaviors  in the case of numerical
 instability, as spurious oscillations may  propagate backward in
 time.

 Our aim in this paper is to design a Spline Upwind (SU) formulation of the heat equation, with 
 a stabilization term that promotes  causality. The proposed  SU 
 generalizes  classical upwinding, as  SUPG (\cite{Brooks1982}),
 to higher degree splines. We recall that SUPG method in time for the heat equation and with piecewise
linear finite elements leads to a lower block-triangular linear system. Stability is further enhanced by adding
artificial diffusion when the residual is higher, as with Shock Capturing. These techniques concomitantly
promotes causality and stability, thereby enhancing the overall computational robustness. The proposed
SU extends these ideas to higher degree splines. We first enrich the plain Galerkin formulation by diffusion
terms of different order, such that the resulting linear system is block triangular with respect to time. These
terms are then weighted by the residual in order to preserve the optimal convergence rate when the solution
is smooth.

We perform numerical tests to assess the expected  behavior of the proposed SU formulation. In particular, motivated
by the interest for space-time simulation of laser-based additive
manufacturing (see \cite{Kopp2021}), we perform experiments with a
concentrated source term, showing that the numerical solution is
free from spurious oscillations.

While the focus of our work is not on
  computational cost, we also acknowledge that space-time
  formulations pose challenges in terms of their computational
  cost.  The augmented dimensionality detrimentally impacts
  conventional solvers. However, it is noteworthy that space-time
  formulations hold promise for local mesh refinement
  \cite{Langer2021}  and parallelisation
  \cite{Gander2015}, attracting 
  interest in the field, see also   the recent book \cite{langer2019}.

The outline of the paper is as follows.  The basics of IgA are discussed in Section \ref{sec:preliminaries}. In Section
\ref{sec:adv_problem} we review some stabilized formulations in one
dimension, for advection and advection-diffusion equations, and
introduce the new SU formulation.
In Section \ref{sec:heat_problem} we apply  SU to the heat equation.
We propose numerical tests, assessing the performance of the presented stabilizing methods, in Section \ref{sec:results}.
Finally, in the last section we draw conclusions and highlight
some future research directions.


\section{Preliminaries} \label{sec:preliminaries}

We recall the notation and definitions of \cite{Loli2020}.

Given $n$ and $p$ two positive integers, we consider the knot vector 
\begin{equation*}
	\widehat{\Xi}:=\left\{ 0=\widehat{\xi}_1=\dots=\widehat{\xi}_{p+1} \leq \dots \leq \widehat{\xi}_{n}=\dots=\widehat{\xi}_{n+p+1}=1\right\}
\end{equation*}
and the vector $\widehat{Z}:=\left\{\widehat{\zeta}_1, \dots, \widehat{\zeta}_m\right\}$ of knots without repetitions (i.e. breakpoints).

The univariate spline space is defined as
\begin{equation*}
	\widehat{\S}_h^p : = \mathrm{span}\{\widehat{b}_{i,p}\}_{i = 1}^n,
\end{equation*}
where $\widehat{b}_{i,p}$ are the univariate B-splines and $ \widehat{h}$ denotes the mesh-size, i.e. $ \widehat{h}:=\max\{ |\widehat{\xi}_{i+1}-\widehat{\xi}_i| \ \text{s.t.} \ i=1,\dots,n+p \}$. 
For more details on B-splines properties  and their use in IgA we refer to \cite{Cottrell2009}.

Multivariate B-splines are tensor product of univariate B-splines.  
We consider  functions that depend on $d$ spatial variables  and the time variable. 
Given positive integers $n_l, p_l$ for $l=1,\dots,d$ and $n_t,p_t$, we define $d+1$ univariate knot vectors $\widehat{\Xi}_l:=\left\{ \widehat{\xi}_{l,1} \leq \dots \leq \widehat{\xi}_{l,n_l+p_l+1}\right\}$  for $l=1,\ldots, d$ and $\widehat{\Xi}_t:=\left\{\widehat{\xi}_{t,1} \leq \dots \leq \widehat{\xi}_{t,n_t+p_t+1}\right\}$ and $d+1$ breakpoints vectors $\widehat{Z}_l:=\left\{ \widehat{\zeta}_{l,1},\dots,\widehat{\zeta}_{l,m_l}\right\}$  for $l=1,\ldots, d$ and $\widehat{Z}_t:=\left\{\widehat{\zeta}_{t,1},\dots,\widehat{\zeta}_{t,m_t}\right\}$.
Let $\widehat{h}_l$ be the mesh-size associated to the knot vector $\widehat{\xi}_l$ for $l=1,\dots,d$, let $\widehat{h}_s:=\max\{\widehat{h}_l\ | \ l=1,\dots, d\}$ be the maximal mesh-size in all spatial knot vectors and let $\widehat{h}_t$ be the mesh-size of the time knot vector.
	
Let also $\boldsymbol{p}$ be the vector that contains the degree of each univariate spline space, i.e. $\boldsymbol{p} :=(\boldsymbol{p}_s,
p_t)$, where $\boldsymbol{p}_s:= (p_1,\dots,p_d )$. 
 
The multivariate B-splines are defined as
\begin{equation*} 
	\widehat{B}_{ \vect{i},\vect{p}}(\vect{\eta},\tau) : = \widehat{B}_{\vect{i_s}, \vect{p}_s}(\vect{\eta}) \widehat{b}_{i_t,p_t}(\tau),
\end{equation*}
where 
\begin{equation*} 
  \widehat{B}_{\vect{i_s},\vect{p}_s}(\vect{\eta}):=\widehat{b}_{i_1,p_1}(\eta_1) \ldots \widehat{b}_{i_d,p_d}(\eta_d),
\end{equation*}
$\vect{i_s}:=(i_1,\dots,i_d)$, $\vect{i}:=(\vect{i_s}, i_t)$  and  $\vect{\eta} = (\eta_1, \ldots, \eta_d)$.  
The corresponding spline space is defined as
\begin{equation*}
	\widehat{\S} ^{\vect{p}}_{ {h}  }  := \mathrm{span}\left\{\widehat{B}_{\vect{i}, \vect{p}} \ \middle| \ i_l = 1,\dots, n_l \text{ for } l=1,\dots,d; i_t=1,\dots	,n_t \right\},
\end{equation*} 
and $\widehat{h}:=\max\{\widehat{h}_s,\widehat{h}_t\}$. 
We have that
$\widehat{\S} ^{\vect{p}}_{ {h}}  =\widehat{\S} ^{ \vect{p}_s}_{ {h}_s} \otimes \widehat{\S} ^{p_t}_{h_t}, $ where   \[\widehat{\S} ^{\vect{p}_s}_{h_s} := \mathrm{span}\left\{\widehat{B}_{\vect{i_s},\vect{p}_s} \ \middle|  \ i_l =  1,\dots, n_l; l=1,\dots,d  \right\}\] is the space of tensor-product splines on $\widehat{\Omega}:=(0,1)^d$.  

We assume that $p_t, p_s\geq 1$ and that $\widehat{\S} ^{\vect{p}_s}_{h_s}
	\subset C^0(\widehat{\Omega}  )$ and  $\widehat{\S} ^{{p}_t}_{h_t}
	\subset C^{{p}_t-1}\left((0,1)\right)$. We allow variable
        continuity in space since it may be useful for  geometry
        representation, while we consider only maximum continuity with
        respect to time in order to benefit from the approximation
        properties of smooth splines, see \cite{Evans2009,Bressan2018}.

We denote by $\Omega\times (0,T)$ the space-time computational domain, where $\Omega\subset\mathbb{R}^d$ ($d$ denotes the space dimension) and $\Omega$ is parametrized by  $\vect{F}: \widehat{\Omega} \rightarrow {\Omega}$, with  $\vect{F}\in  {\widehat{\mathcal{S}}^{\vect{p}_s}_{{h}_s}}$, and $T>0$ is the final time.
The space-time domain is parametrized by $\vect{G}:\widehat{\Omega}\times(0,1)\rightarrow
 \Omega\times(0,T)$, such that $ \vect{G}(\vect{\eta}, \tau):=(\vect{F}(\vect{\eta}), T\tau )=(\vecx,t).$
 
The spline space with initial and boundary conditions, in parametric coordinates, is  
\begin{equation*}
\widehat{\mathcal{X}}_{h}:=\left\{ \widehat{v}_h\in \widehat{\mathcal{S}}^{\vect{p}}_h \ \middle| \ \widehat{v}_h = 0 \text{ on } \partial\widehat{\Omega}\times (0,1) \text{ and } \widehat{v}_h = 0 \text{ on } \widehat{\Omega}\times\{0\} \right\}.
\end{equation*}
  We also have that  $
\widehat{\mathcal{X}}_{h} =   \widehat{\mathcal{X}}_{s,h_s}     \otimes  \widehat{\mathcal{X}}_{t,h_t}  $,  where 
 \begin{align*}
  \widehat{\mathcal{X}}_{s,h_s}   & := \left\{ \widehat{w}_h\in \widehat{\mathcal{S}}^{\vect{p}_s}_{h_s}    \ \middle| \ \widehat{w}_h = 0 \text{ on } \partial\widehat{\Omega}  \right\}\  \\
  &  \ = \ \text{span}\left\{ \widehat{b}_{i_1,p_s}\dots\widehat{b}_{i_d,p_s} \ \middle| \ i_l = 2,\dots , n_l-1; \ l=1,\dots,d\ \right\},\\ 
   \widehat{\mathcal{X}}_{t,h_t} & := \left\{ \widehat{w}_h\in \widehat{\mathcal{S}}^{ p_t}_{h_t} \ \middle|  \ \widehat{w}_h( 0)=0 \right\}  \   = \ \text{span}\left\{ \widehat{b}_{i_t,p_t} \ \middle| \ i_t = 2,\dots , n_t\ \right\}.
 \end{align*}
With a colexicographical reordering of the basis functions, we write  
\begin{align*}
  \widehat{\mathcal{X}}_{s,h_s} &   = \ \text{span}\left\{ \widehat{b}_{i_1,p_s}\dots\widehat{b}_{i_d,p_s} \ \middle| \ i_l = 1,\dots , n_{s,l}; \ l=1,\dots,d\ \right\}\\ 
  & \ =\text{span}\left\{ \widehat{B}_{i, \vect{p}_s} \ \middle|\ i =1,\dots , N_s   \ \right\},\\
\  \widehat{\mathcal{X}}_{t,h_t}  &  = \ \text{span}\left\{ \widehat{b}_{i,p_t} \ \middle| \ i = 1,\dots , N_t\ \right\},
\end{align*}
  and
\begin{equation}
\widehat{\mathcal{X}}_{h}=\text{span}
\left\{ \widehat{B}_{{i}, \vect{p}} \ \middle|\ i=1,\dots,N_{dof} \right\},
\label{eq:all_basis}
\end{equation}
where $n_{s,l}:= n_l-2 $ for $l=1,\dots,d$, $N_s:=\prod_{l=1}^dn_{s,l}$, $N_t:=n_t-1$ and $ N_{dof}:=N_s N_t$. 

Our isogeometric space is the isoparametric push-forward of \eqref{eq:all_basis} through the geometric map $\vect{G}$, i.e.
\begin{equation*}
\mathcal{X}_{h} := \text{span}\left\{  B_{i, \vect{p}}:=\widehat{B}_{i, \vect{p}}\circ \vect{G}^{-1} \ \middle| \ i=1,\dots , N_{dof}   \right\},
\end{equation*}
where again  
 $\mathcal{X}_{h}=\mathcal{X}_{s,h_s}   \otimes \mathcal{X}_{t,h_t} $,
   with 
\begin{equation*} 
 \mathcal{X}_{s,h_s}    :=\text{span}\left\{ {B}_{i, \vect{p}_s}:= \widehat{B}_{i, \vect{p}_s}\circ \vect{F}^{-1} \ \middle| \ i=1,\dots,N_s \right\}
\end{equation*}
and 
\begin{equation*} \label{eq:time_disc_space}
 \mathcal{X}_{t,h_t}   :=\text{span}\left\{  {b}_{i,p_t}(\cdot):= \widehat{b}_{i,p_t}\left(\frac{\cdot}{T}\right) \ \middle| \ i=1,\dots,N_t \right\}.
\end{equation*}

Moreover we define the breakpoints in the time
  interval as: 
\begin{equation*} 
 \zeta_{t,i}:= T\widehat{\zeta}_{t,i} \ \ \ \text{for} \ \ i=1,\cdots,m_t,
\end{equation*} and the time steps as: 
\begin{equation*} 
 h_{t,i}:= \zeta_{t,i+1}-\zeta_{t,i}\ \ \ \text{for} \ \ i=1,\cdots,m_t-1.
\end{equation*}


\section{Upwinding in one dimension} \label{sec:adv_problem}

Our first model problem is the unidimensional advection: we look for a function $u$ such that 
\begin{equation}
\label{eq:adv_problem}
	\left\{
	\begin{array}{rlcll}
		  &  u' & = & f & \mbox{in }  (0, T)\\[1pt]
			&  u(0) & = & 0 &
	\end{array}
	\right.
\end{equation}
We assume that  $f\in L^2(0,T)$ and consider the following Galerkin method: 
\begin{equation*}
 \text{find }   u_h\in \mathcal{X}_{t,h_t} \text{ such that } \mathcal{A}(u_h; v_h) = \mathcal{F}(v_h) \quad \, \forall v_h\in  \mathcal{X}_{t,h_t},
\end{equation*} 
where
\[
 \mathcal{A}(u_h;v_h) := \int_{0}^T   u'_h \,  v_h	\, \dt  
\quad \text{and} \quad
\mathcal{F}(v)  := \int_0^T f\, v_h \, \dt.
\] 


\subsection{Standard Upwind and Shock Capturing} \label{sec:standard_stab}

The Streamline Upwind Petrov Galerkin (SUPG) method, see \cite{Brooks1982}, reads:
\begin{equation}
\label{eq:adv_SUPG_discrete-system}
 \text{find }   u_h\in \mathcal{X}_{t,h_t} \text{ such that } \mathcal{A}(u_h; v_h) + \mathcal{S}_{\text{SUPG}}(u_h, f; v_h)= \mathcal{F}(v_h) \quad \, \forall v_h\in  \mathcal{X}_{t,h_t},
\end{equation} 
where 
\begin{equation*}
 \mathcal{S}_{\text{SUPG}}(u_{h}, f; v_h) :=  \sum_{i=1}^{m_t-1}\tau_{\text{SUPG},i} \int_{\zeta_{t,i}}^{\zeta_{t,i+1}}{( u'_{h}-f)   v'_{h}} \ \dt, 
\end{equation*}

In order to enhance the stability of SUPG, we can further add  a
Shock Capturing term  leading to:
\begin{equation*}
 \text{find }   u_h \in \mathcal{X}_{t,h_t} \text{ such that } \mathcal{A}(u_h ; v_h) + \mathcal{S}_{\text{SUPG}}(u_h , f; v_h) + \mathcal{S}_{\text{SC}} (u_h ;v_h) = \mathcal{F}(v_h) \quad \, \forall v_h\in  \mathcal{X}_{t,h_t},
\end{equation*} 
where, following \cite{Bazilevs2007},
\begin{equation*}
 \mathcal{S}_{\text{SC}} (u_h ;v_h) := \sum_{i=1}^{m_t-1}\int_{\zeta_{t,i}}^{\zeta_{t,i+1}} \kappa_{\text{SC},i}  u'_{h} v'_{h} \ \dt \quad \, \text{and} \ \quad \kappa_{\text{SC},i} :=  \tau_{\text{SC},i}\frac{\left| { u'_{h}} - f \right|}{u_{\text{ref}}}, \ \ \ \text{with} \ \tau_{\text{SC},i}:=\frac{h_{t,i}^2}{4} 
\end{equation*} 
and $u_{\text{ref}}$ is a reference magnitude for $u_{h}$.


\subsection{High-order Upwind} \label{sec:Upwind}
\begin{figure}[htbp]
	\centering
	\begin{subfigure}{0.75\textwidth}
			\includegraphics[width=1.00\textwidth]{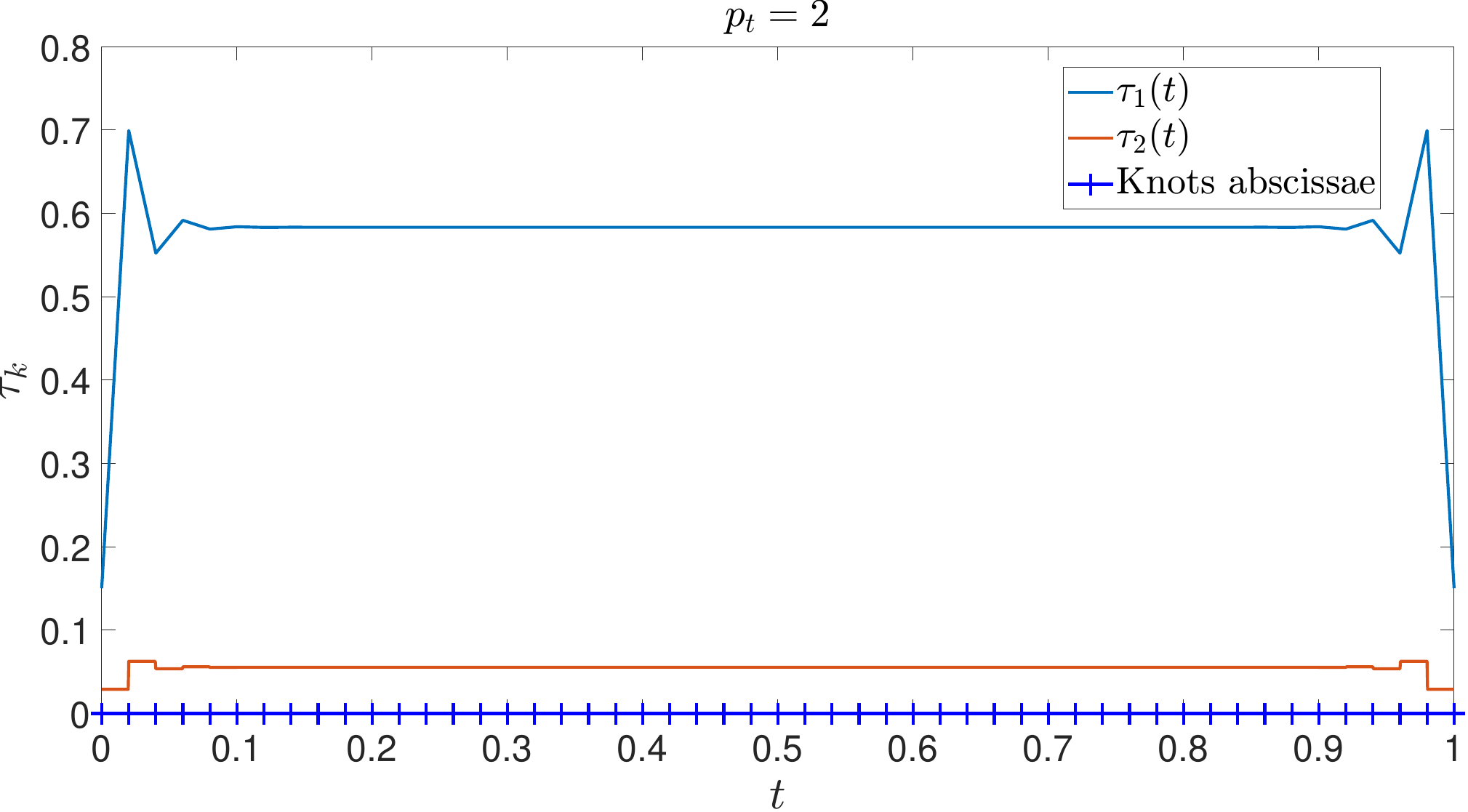}
		\end{subfigure}
	\begin{subfigure}{0.75\textwidth}
			\includegraphics[width=1.00\textwidth]{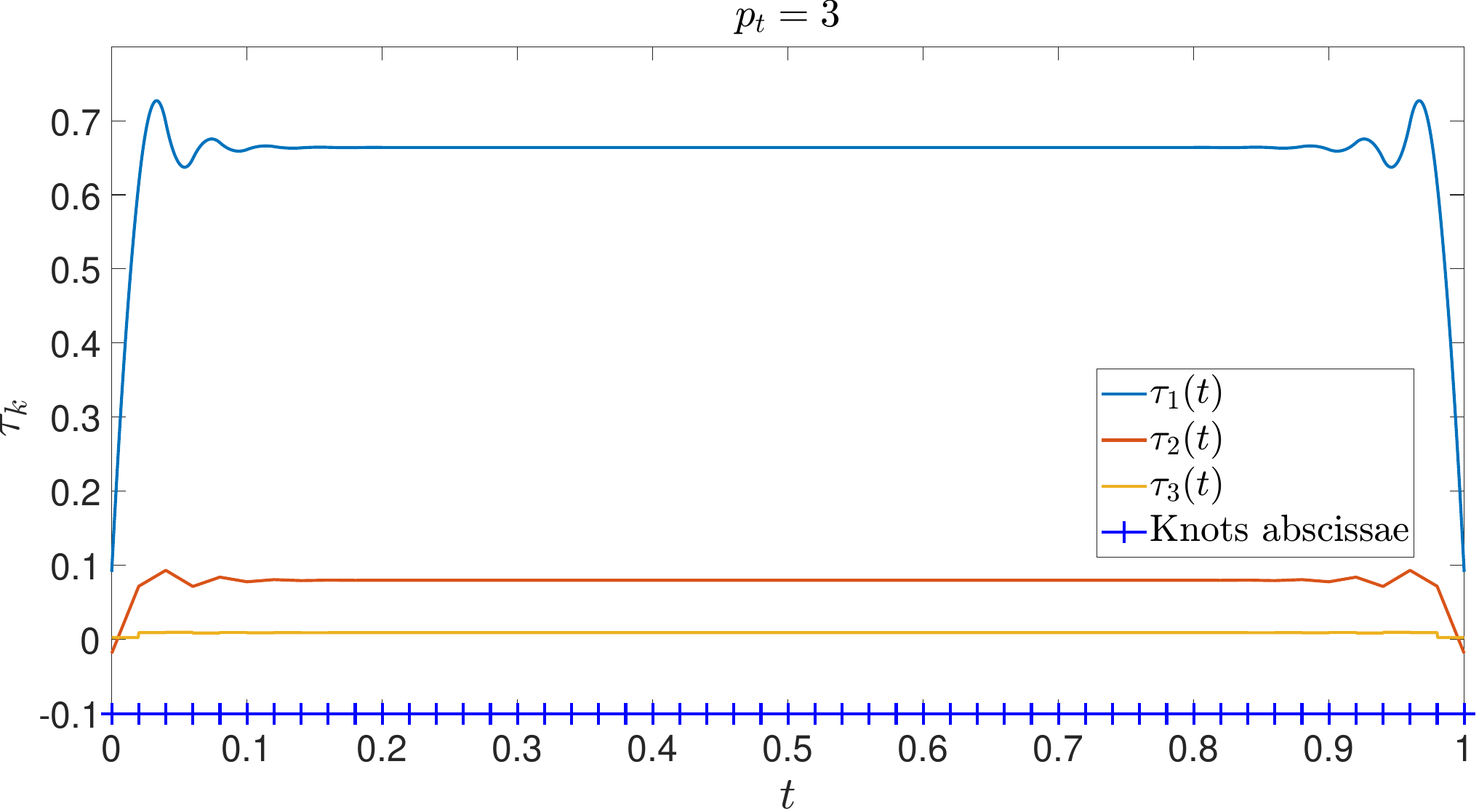}
	\end{subfigure}
	\begin{subfigure}{0.75\textwidth}
			\includegraphics[width=1.00\textwidth]{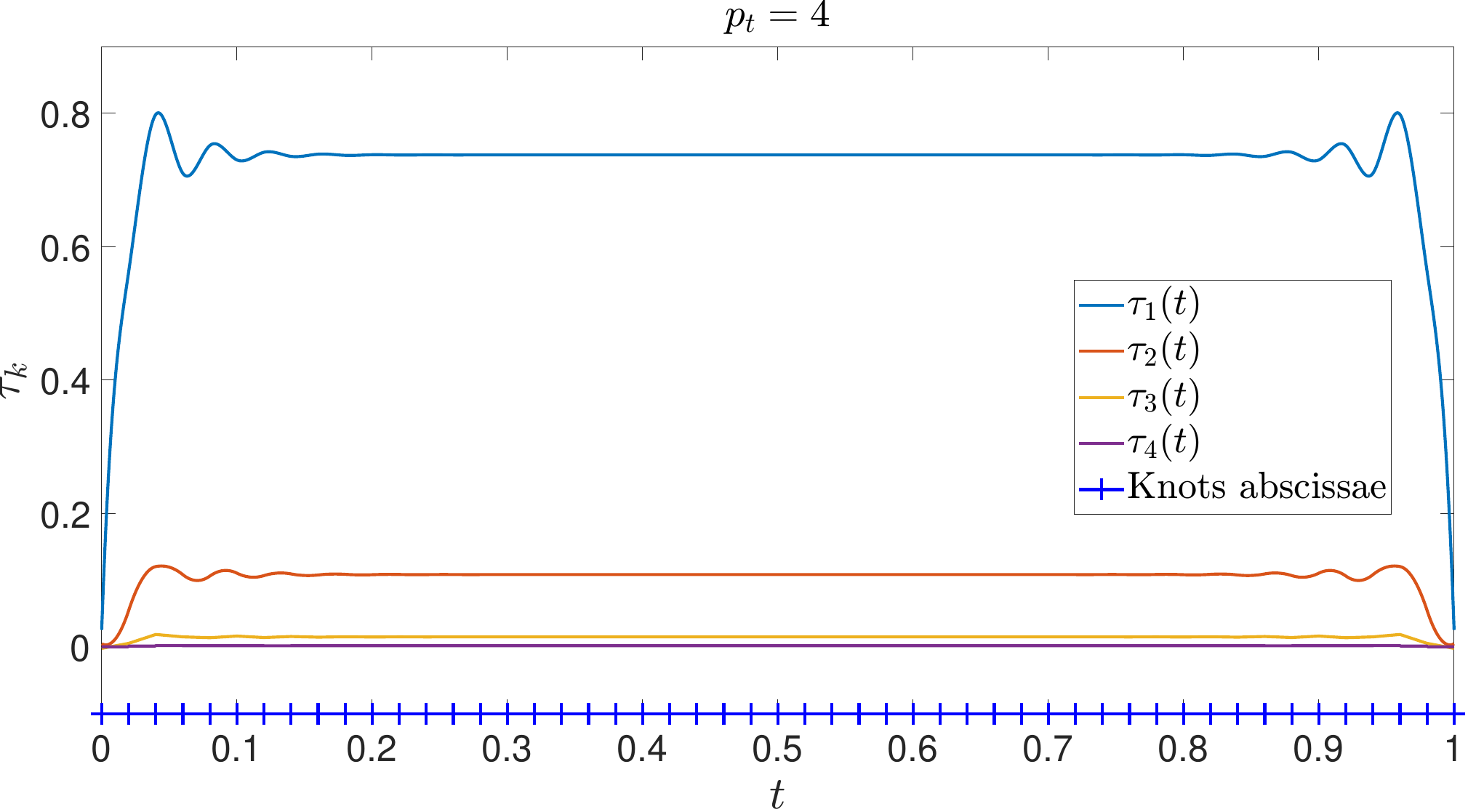}
	\end{subfigure}
	\caption{Plot of $\tau_{k}(t)$ for uniform meshes,
          with $h_t=1/50$ and $p_t=2,3,4$.}
	\label{fig:tau_unif}
\end{figure}

\begin{figure}[htbp]
	\centering
		\includegraphics[width=0.95\textwidth]{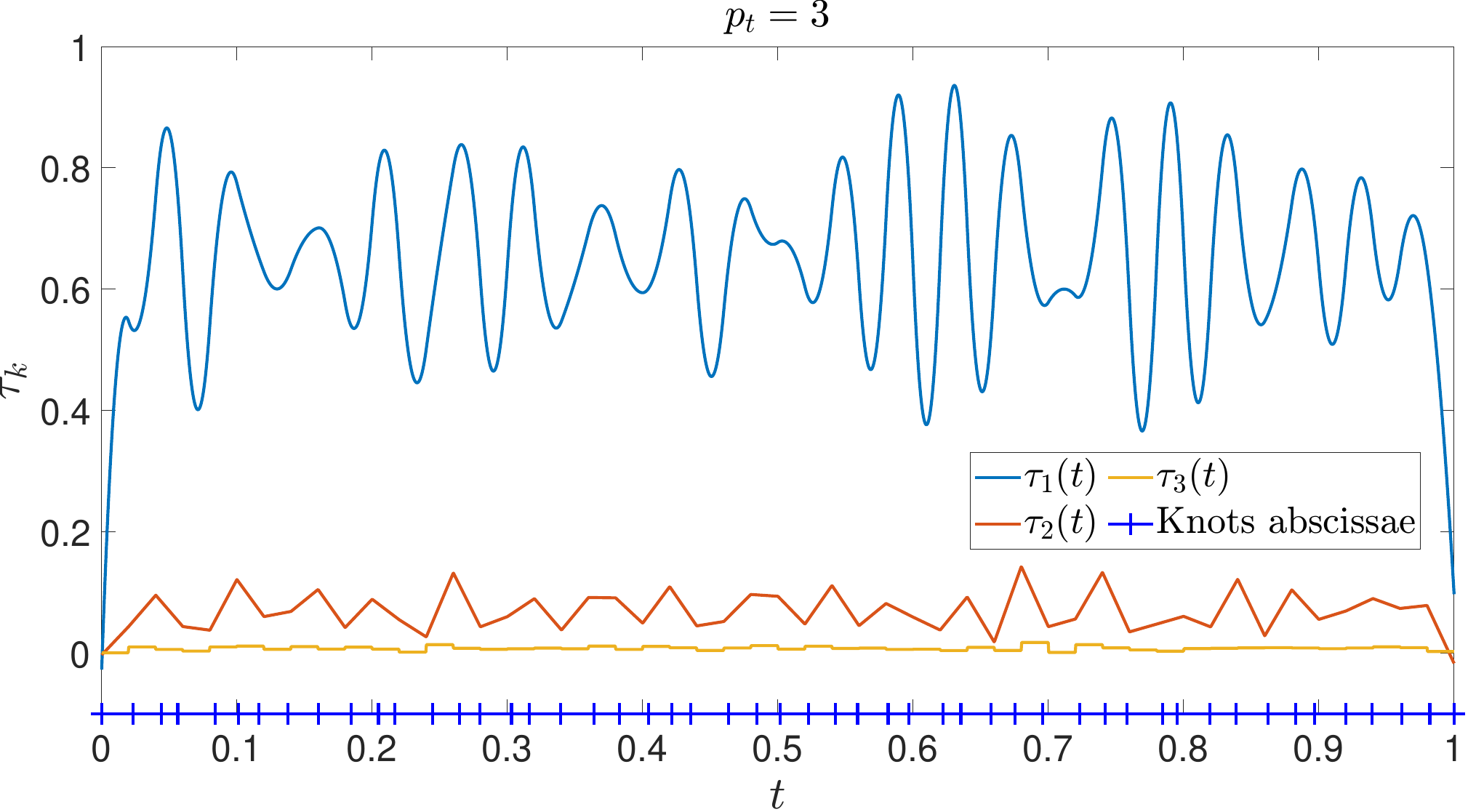}
	\caption{Plot of $\tau_{k}(t)$ for a non-uniform mesh
          (the knots are represented by the blue vertical segments on
          the horizontal axis) and $p_t=3$.}
	\label{fig:tau_nunif}
\end{figure}

With the choice $\tau_{\text{SUPG},i} = \frac{h_{t,i}}{2}$ and when $p_t = 1$, formulation \eqref{eq:adv_SUPG_discrete-system} leads to a lower triangular linear system. However, for higher degree splines, the matrix does not exhibit a lower triangular structure regardless of the $\tau_{\text{SUPG},i}$ value chosen.
This motivates the design of  a new high-order Non-Consistent Spline Upwind (NCSU) formulation for spline with maximum continuity $C^{p_t-1}$:

\begin{equation}
\label{eq:adv_ncons_discrete-system}
 \text{find } u_h\in \mathcal{X}_{t,h_t} \text{ such that } \mathcal{A}(u_h; v_h) + \mathcal{S}_{\text{NCSU}}(u_{h}; v_{h})= \mathcal{F}(v_h) \quad \, \forall v_h\in  \mathcal{X}_{t,h_t},
\end{equation}
where the new stabilizing term fulfils
\begin{equation*}
 \mathcal{S}_{\text{NCSU}}(u_{h}; v_{h}) := \sum_{k=1}^{p_t}\sum_{i=1}^{m_t-1}{h_{t,i}}^{2k-1} \int_{\zeta_{t,i}}^{\zeta_{t,i+1}}{\tau_k(t){u_{h}^{(k)} v_{h}^{(k)}} \ \dt},
\end{equation*}
where each $\tau_{k}(T \ \cdot) \in \widehat{\mathcal{S}}^{
  p_t-k}_{h_t} \in C^{p_t-k-1}$ is a spline with maximum continuity
that is selected in order to make the linear  system matrix lower
triangular, that is:
\begin{equation}\label{sys_tau}
	\int_0^T{ b'_{\ell+i,p_{t}} b_{i,p_{t}} \ \dt} + \sum_{k=1}^{p_t}\sum_{j=1}^{m_t-1}{h_{t,j}}^{2k-1} \int_{\zeta_{t,j}}^{\zeta_{t,j+1}}{\tau_k(t){ b^{(k)}_{\ell+i,p_{t}} b^{(k)}_{i,p_{t}}} \ \dt} = 0,
\end{equation}
for $i = 1,\ldots,N_t-1\ \text{and} \ \ \ell=1,\ldots,r, \  \text{with} \ \ r=\min(p_t,N_t-i)$.\\
In all our numerical tests we have observed that \eqref{sys_tau} is
well posed and that  the $\tau_{k}$ are bounded; in particular the
$\tau_{k}$ are positive in most of the domain and where they are
negative, they are much smaller in magnitude.  In Figures \ref{fig:tau_unif} and \ref{fig:tau_nunif} we have plotted $\tau_k(t)$ for different degrees and meshes.

However, the formulation above is non consistent and lacks optimal convergence.
To overcome this limitation we introduce a non-linear weighting based on residual
denoted for brevity Spline Upwind (SU) method, and
defined as:
\begin{equation}
\label{eq:adv_uim_discrete-system}
 \text{find }   u_h \in \mathcal{X}_{t,h_t} \text{ such that } \mathcal{A}(u_h ; v_h) + {\mathcal{S}_{\text{SU},1}}(u_h , f; v_h) + {\mathcal{S}_{\text{SU,2}} }(u_h ; v_h) = \mathcal{F}(v_h) \quad \, \forall v_h\in  \mathcal{X}_{t,h_t},
\end{equation} 
with 
\begin{equation*}
	 \mathcal{S}_{\text{SU,1}}(u_h, f; v_h) := \sum_{i=1}^{m_t-1} h_{t,i} \int_{\zeta_{t,i}}^{\zeta_{t,i+1}}{\tau_{1}(t) ( u_h' - (1-\theta (t))f) v_h'} \ \dt,	
\end{equation*}
and
\begin{equation*}
 {\mathcal{S}_{\text{SU,2}} }(u_h ; v_h) := \sum_{k=2}^{p_t} \sum_{i=1}^{m_t-1} { h_{t,i}^{2k-1} \int_{\zeta_{t,i}}^{\zeta_{t,i+1}}{\tau_{k}(t)\theta (t) u_{h}^{(k)} v_{h}^{(k)}} \ \dt},
\end{equation*}
where $\theta (t)$ is a piecewise linear interpolation of $\theta_i $ computed in the breakpoints $\zeta_{t,i}$ for $i=1,\ldots,m_t$ as:
\begin{equation*} 
	\theta_i :=\min (\text{res}_i ,1),
\end{equation*} 
with the relative residual $\text{res}_i$ defined as
\begin{equation*} 
	\text{res}_i :=\frac{\left\| { u'_{h}} - f \right\|_{L^{\infty}([\zeta_{t,\max(1,i-1)},\zeta_{t,\min(i+1,m_t)}])}}{T^{-1}\left\| {u_{h}} \right\|_{L^{\infty}([0,T])} + \left\| {u'_{h}} \right\|_{L^{\infty}([0,T])}}.
\end{equation*} 
The definition of $\theta$ is such that when the residual is high, e.g., within layers, $\theta \approx$ 1 and SU locally reduces to NCSU \eqref{eq:adv_ncons_discrete-system}.

For advection-diffusion problem

\begin{equation}
\label{eq:adv_diff_problem}
	\left\{
	\begin{array}{rcllrcl}
			-\varepsilon  u'' +  u'  & = & f & \mbox{in }  (0, T)\\[1pt]
			 u(0)=u(T) & = & 0 & 
	\end{array}
	\right.
\end{equation}
where $\varepsilon>0$, SU is extended straightforwardly by redefining $\mathcal{S}_{\text{SU,1}}$ in \eqref{eq:adv_uim_discrete-system} as follows
\begin{equation*}
	 \mathcal{S}_{\text{SU,1}}(u_h, f; v_h) := \sum_{i=1}^{m_t-1} h_{t,i} \int_{\zeta_{t,i}}^{\zeta_{t,i+1}}{\tau_{1}(t) ( u_h' - (1-\theta (t))(-\varepsilon u_h'' - f)) v_h'} \ \dt,	
\end{equation*}
and 
\begin{equation*} 
	\theta _{i}:=\min (\text{res}_i,1),
\end{equation*} 
with
\begin{equation*} 
	\text{res} _i:=\frac{\left\| - \varepsilon u''_{h} +  u'_{h} - f \right\|_{L^{\infty}([\zeta_{t,\max(1,i-1)},\zeta_{t,\min(i+1,m_t)}])}}{T^{-1}\left\| {u_{h}} \right\|_{L^{\infty}([0,T])} + \left\|u'_{h} \right\|_{L^{\infty}([0,T])}}.
\end{equation*}


\section{Upwinding the heat equation} \label{sec:heat_problem}

Consider the heat equation with homogeneous boundary and initial condition
\begin{equation}
\label{eq:heat_problem}
	\left\{
	\begin{array}{rcllrcl}
		 \partial_t u -  \Delta u  & = & f & \mbox{in } &\Omega \!\!\!\! &\times & \!\!\!\! (0, T) \\[1pt]
		 u  & = & 0 & \mbox{on } &\partial\Omega \!\!\!\! &\times& \!\!\!\! [0, T] \\[1pt]
		 u & = & 0 & \mbox{in } &\Omega\!\!\!\!  &\times &\!\!\!\! \lbrace 0 \rbrace
	\end{array}
	\right.
\end{equation}

Introducing the bilinear form $\mathcal{A}(\cdot;\cdot)$ and the linear form $\mathcal{F}(\cdot)$ as
\[
 \mathcal{A}(w;v) := \int_{0}^T\int_{\Omega}  \left( \partial_t {w} \, v +  \nabla w\cdot \nabla v \right)\,\d\Omega\,   \dt  
\quad \text{and} \quad
\mathcal{F}(v)  := \int_0^T\int_{\Omega}f\,   v \,\d\Omega\,\dt
\] 
we consider the Galerkin method: 
\begin{equation*}
 \text{find }   u_h\in \mathcal{X}_h \text{ such that } \mathcal{A}(u_h; v_h) = \mathcal{F}(v_h) \quad \, \forall v_h\in  \mathcal{X}_h.
\end{equation*} 

The matrix of the linear system is
\begin{equation} 
\mathbf{W}_t \otimes  \mathbf{M}_s + \mathbf{M}_t\otimes  \mathbf{K}_s, 
\label{eq:heat_syst_mat} 
\end{equation}
where  for $i,j=1,\dots,N_t$
\begin{subequations}
\begin{equation*}
\label{eq:heat_time_mat}
 [ \mathbf{W}_t]_{i,j}  = \int_{0}^T   b'_{j,  {p}_t}(t)\,  b_{i,{p}_t}(t) \, \dt \quad \text{and} \quad   [ \mathbf{M}_t]_{i,j} = \int_{0}^T\,  b_{j, p_t}(t)\,  b_{i, p_t}(t)  \, \dt , 
  \end{equation*}
  while for $  i,j=1,\dots,N_s $
  \begin{equation*}
  \label{eq:heat_space_mat}
 [ \mathbf{K}_s]_{i,j}  =  \int_{\Omega} \nabla  B_{j, \vect{p}_s}(\vect{x})\cdot \nabla  B_{i, \vect{p}_s}(\vect{x}) \ \d\Omega  \quad \text{and} \quad  [ \mathbf{M}_s]_{i,j}  =  \int_{\Omega}  B_{j, \vect{p}_s}(\vect{x}) \  B_{i, \vect{p}_s}(\vect{x}) \ \d\Omega .  
 \end{equation*}
 \label{eq:pencils}
\end{subequations}

The SUPG method reads:
\begin{equation}
\label{eq:heat_SUPG_discrete-system}
 \text{find }   u_h\in \mathcal{X}_h \text{ such that } \mathcal{A}(u_h; v_h) + {\mathcal{S}_{\text{SUPG}}}(u_h, f; v_h)= \mathcal{F}(v_h) \quad \, \forall v_h\in  \mathcal{X}_{h},
\end{equation} 
where
\begin{equation*}
 {\mathcal{S}_{\text{SUPG}}}(u_{h}, f; v_h):= \sum_{i=1}^{m_t-1}\tau_{\text{SUPG},i}\int_{\zeta_{t,i}}^{\zeta_{t,i+1}}\int_{\Omega}{( \partial_t u_h -  \Delta u_h- f)}\partial_t v_{h} \,\d\Omega\,\dt.\\ 
\end{equation*}
With the choice $\tau_{\text{SUPG},i} = \frac{h_{t,i}}{2}$ and when $p_t=1$, formulation \eqref{eq:heat_SUPG_discrete-system} lead to a lower block triangular time derivative matrix.

The new space-time formulation for the heat equation is based, as in the SU one-dimensional formulation, on the idea of modifying \eqref{eq:heat_syst_mat} in order to obtain lower triangular time matrices. 
This is accomplished by introducing numerical
  diffusion in the time direction. The proposed method reads:
\begin{equation}
\label{eq:heat_uim_discrete-system}
 \text{find }   u_h \in \mathcal{X}_{h} \text{ such that } \mathcal{A}(u_h ; v_h) + {{\mathcal{S}_{\text{SU,1}} }}(u_h , f; v_h) + {{\mathcal{S}_{\text{SU,2}} }}(u_h ; v_h) + {{\mathcal{S}_{\text{SU,3}} }}(u_h ;v_h)= \mathcal{F}(v_h) \quad \, \forall v_h\in  \mathcal{X}_{h},
\end{equation} 
where for $i,j = 1,\ldots,N_{dof}$
\begin{equation*}
	{{\mathcal{S}_{\text{SU,1}} }}(u_h , f; v_h) :=  \sum_{i=1}^{m_t-1}h_{t,i} \int_{\zeta_{t,i}}^{\zeta_{t,i+1}} \tau_{1}(t) \int_{\Omega} {(\partial_t u_h + (1-\theta ({\vect{x}},t)) ( -\Delta u_h -  f)) \partial_t v_h} \ \d\Omega \ \dt,
\end{equation*}
and
\begin{equation*}
 {{\mathcal{S}_{\text{SU,2}} }}(u_h ; v_h) := \sum_{k=2}^{p_t}\sum_{i=1}^{m_t-1}{h_{t,i}^{2k-1} \int_{\zeta_{t,i}}^{\zeta_{t,i+1}} \tau_{k}(t) \int_{\Omega} {\theta (\vect{x},t)\partial_t^k u_h\partial_t^k v_h} \ \d\Omega \ \dt},
\end{equation*}
while
\begin{equation*}
 {{\mathcal{S}_{\text{SU,3}} }}(u_h ; v_h) := \sum_{k=1}^{p_t}\sum_{i=1}^{m_t-1}{h_{t,i}^{2k} \int_{\zeta_{t,i}}^{\zeta_{t,i+1}} \sigma_k(t)\int_{\Omega} {\theta (\vect{x},t)(\nabla(\partial_t^k u_h)\cdot\nabla(\partial_t^k v_h})) \ \d\Omega \ \dt}.
\end{equation*}
As in the Section \ref{sec:Upwind}, $\tau_{k}(T \ \cdot) \in \widehat{\mathcal{S}}^{ p_t-k}_{h_t}$ with maximum continuity, are selected such that:
\begin{equation}
	\int_0^T{ b'_{\ell+i,p_{t}} b_{i,p_{t}} \ \dt} + \sum_{k=1}^{p_t}\sum_{j=1}^{m_t-1}{h_{t,j}}^{2k-1} \int_{\zeta_{t,j}}^{\zeta_{t,j+1}}{\tau_k(t){ b^{(k)}_{\ell+i,p_{t}} b^{(k)}_{i,p_{t}}} \ \dt} = 0,
\end{equation}
for $i = 1,\ldots,N_t-1\ \text{and} \ \ \ell=1,\ldots,r, \  \text{with} \ \ r=\min(p_t,N_t-i)$,\\
while $\sigma_{k}(T \ \cdot) \in \widehat{\mathcal{S}}^{ p_t-k}_{h_t}$ with maximum continuity, are selected in order to make the time mass matrix lower triangular:
\begin{equation}
	\int_0^T{ b_{\ell+i,p_{t}} b_{i,p_{t}} \ \dt} + \sum_{k=1}^{p_t}\sum_{j=1}^{m_t-1}{h_{t,j}}^{2k} \int_{\zeta_{t,j}}^{\zeta_{t,j+1}}{\sigma_k(t){ b^{(k)}_{\ell+i,p_{t}} b^{(k)}_{i,p_{t}}} \ \dt} = 0,
\end{equation}
for $i = 1,\ldots,N_t-1\ \text{and} \ \ \ell=1,\ldots,r, \  \text{with} \ \ r=\min(p_t,N_t-i)$.

The function $\theta (\vect{x},t)$ ranges from $0$ to $1$. If we set
$\theta$ as a fixed parameter equal to $1$,
\eqref{eq:heat_uim_discrete-system} yields a block lower triangular
global system matrix. However, in order to achieve optimal order of
convergence, similar to the one-dimensional case (see Section
\ref{sec:Upwind}), we define $\theta (\vect{x},t)$ as a  piecewise
$(d+1)$-linear interpolation of $\theta_{\boldsymbol{i},j} $ computed
in the breakpoints, where  
for $\boldsymbol{i}=(i_1,\dots,i_{d})$, $i_l=1,\dots,m_l$, $l=1,\dots,d$, and $j=1,\dots,m_t$ we set
\begin{equation*}
{\theta }_{\boldsymbol{i},j}:=\min (\text{res}_{\boldsymbol{i},j},1),
\end{equation*}
with
\begin{equation*}
\text{res} _{\boldsymbol{i},j}:=\frac{\left\| \partial_t {u_{h}} - \Delta u_{h} - f \right\|_{L^{\infty}(\psi_s\times\psi_t)}}{T^{-1}\left\| {u_{h}} \right\|_{L^{\infty}(\Omega \times[0,T])} + \left\| \partial_t {u_{h}} \right\|_{L^{\infty}(\Omega \times[0,T])}},
\end{equation*}
where
\begin{equation*}
\psi_s=[\zeta_{1,\max(1,i_1-1)},\zeta_{1,\min(i_1+1,m_1)}]\times\ldots\times [\zeta_{d,\max(1,i_d-1)},\zeta_{d,\min(i_d+1,m_d)}],
\end{equation*}
and
\begin{equation*}
\psi_t=[\zeta_{t,\max(1,j-1)},\zeta_{t,\min(j+1,m_t)}].
\end{equation*}


\section{Numerical Results} \label{sec:results}

In the following, all numerical tests are conducted using Matlab R2023a and the GeoPDEs toolbox \cite{Vazquez2016}.
Just for the sake of simplicity,  in all our tests we consider splines
of the same polynomial degree in all parametric directions for space
and time. Specifically, we set $p_1=\dots=p_d=p_t:=p$. Additionally,
although the proposed methods are designed for maximum regularity only
with respect to time, we choose to use splines of global maximum
continuity $C^{p-1}$ also with respect to space. Numerical tests with
different degrees and regularities (in space) indeed yield results entirely analogous to those
reported below.

Nonlinearities in the equations are addressed through fixed point iterations, and the resulting linear systems are solved using the direct solver provided by Matlab.

We would like to emphasize that the main focus of this work is not on computational costs, and as such, we do not discuss or analyze the efficiency associated with the proposed method.


\subsection{Advection equation} \label{sec:adv_results}

We consider the advection equation \eqref{eq:adv_problem} on $(0,T)$ with $T=1$ and uniform mesh.

\subsubsection{Smooth solution}

\begin{figure}[htbp]
	\centering
		\includegraphics[width=1.00\textwidth]{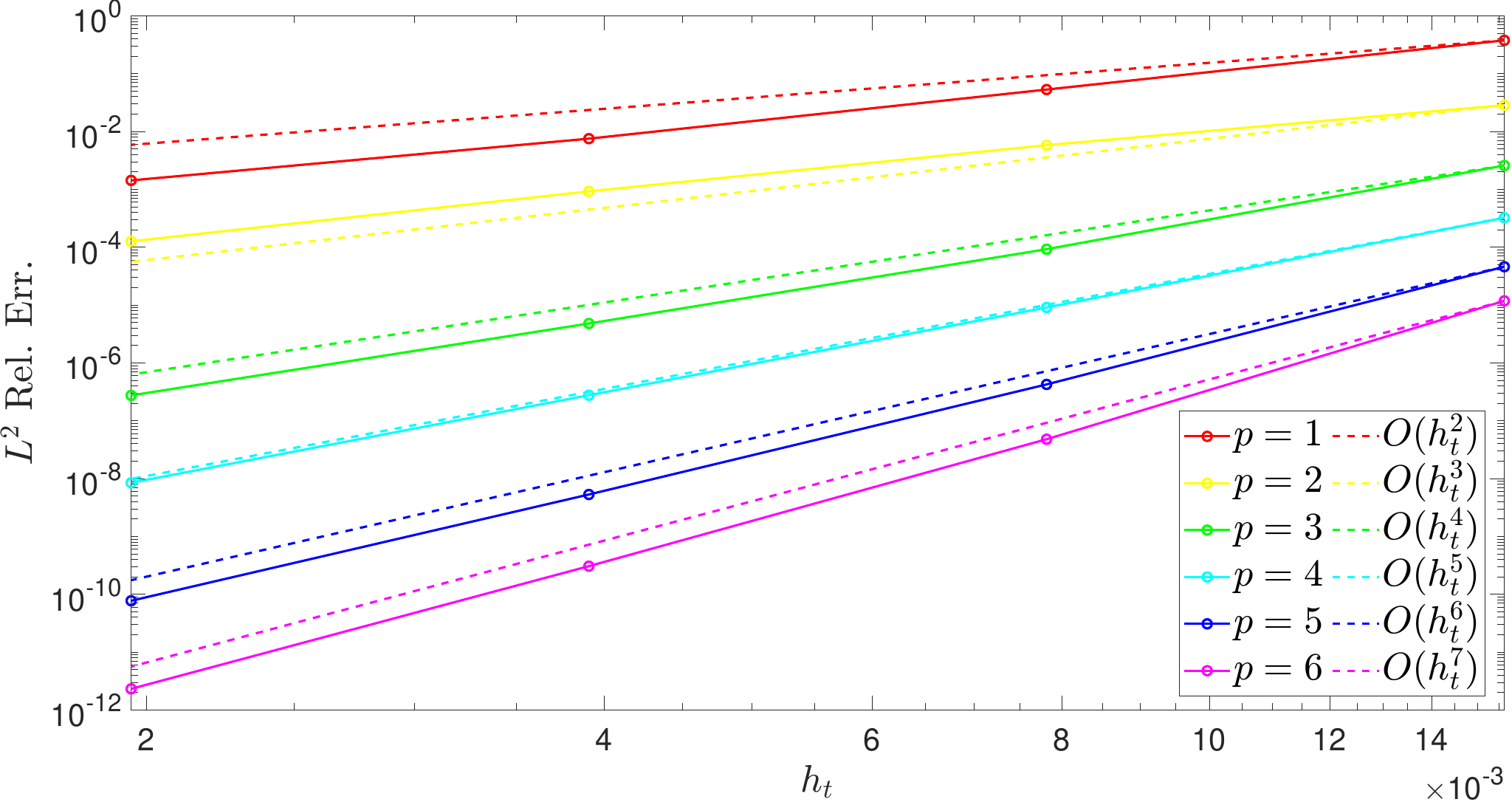}
	\caption{Advection equation, SU relative error plots in $L^2$-norm.}
	\label{fig:SU_L2_error_smooth}
\end{figure}

We set $f=50\cos(50t)$ such that the exact solution is $u_{\text{ex}}(t)=\sin(50t)$.
In Figure \ref{fig:SU_L2_error_smooth}, we show the error plot for the SU formulation on uniform meshes and degree $p=1,\ldots,6$ and we see that the method is optimally convergent.

\subsubsection{Solution with layers}

\begin{figure}[htbp]
	\centering
		\includegraphics[width=1.00\textwidth]{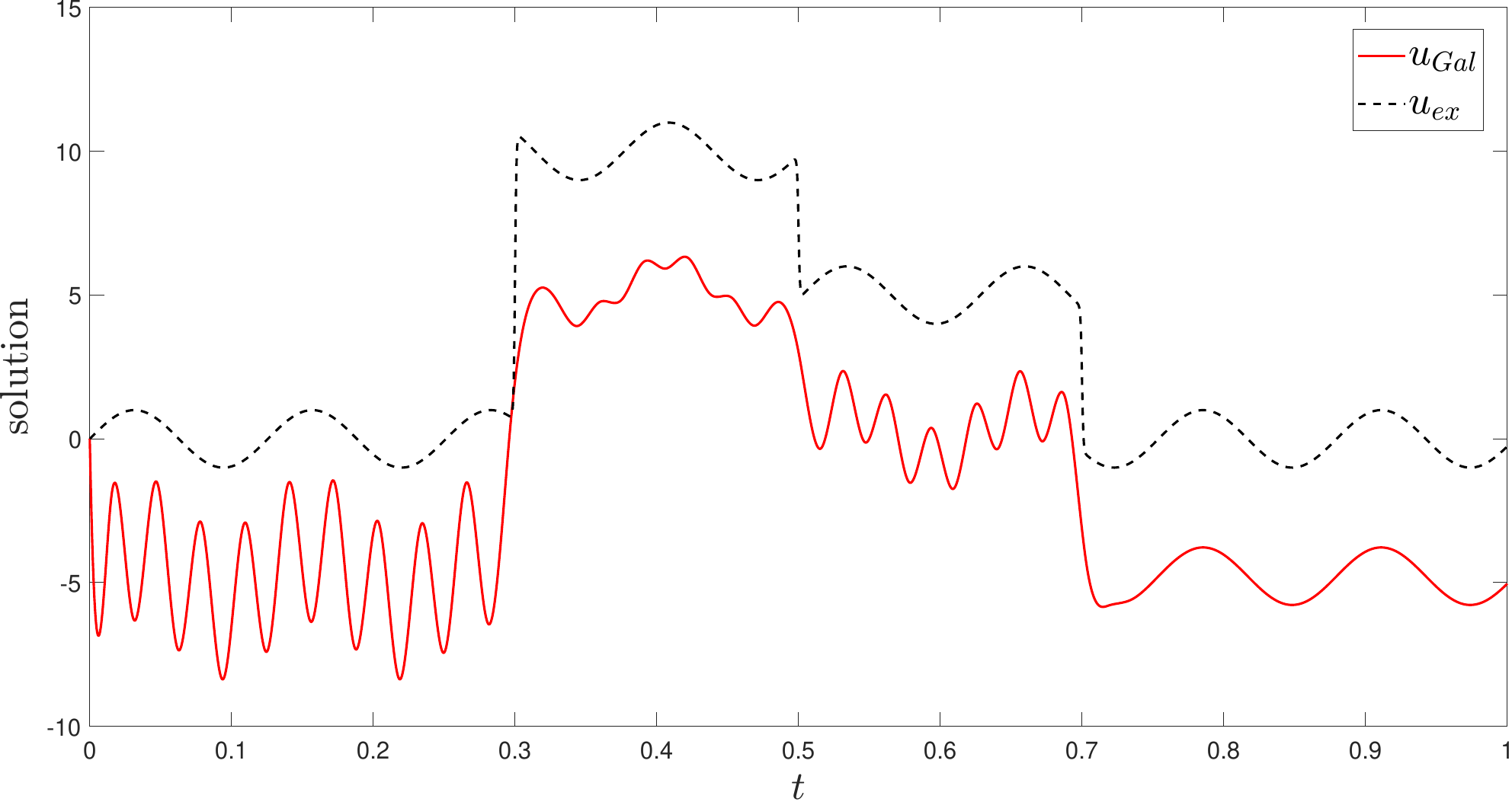}
	\caption{Advection equation, exact and standard Galerkin solutions, with $h_{t}=2^{-6}$ and $p=3$.}
	\label{fig:GALERKIN_layers}
\end{figure}

\begin{figure}[htbp]
	\centering
		\includegraphics[width=1.00\textwidth]{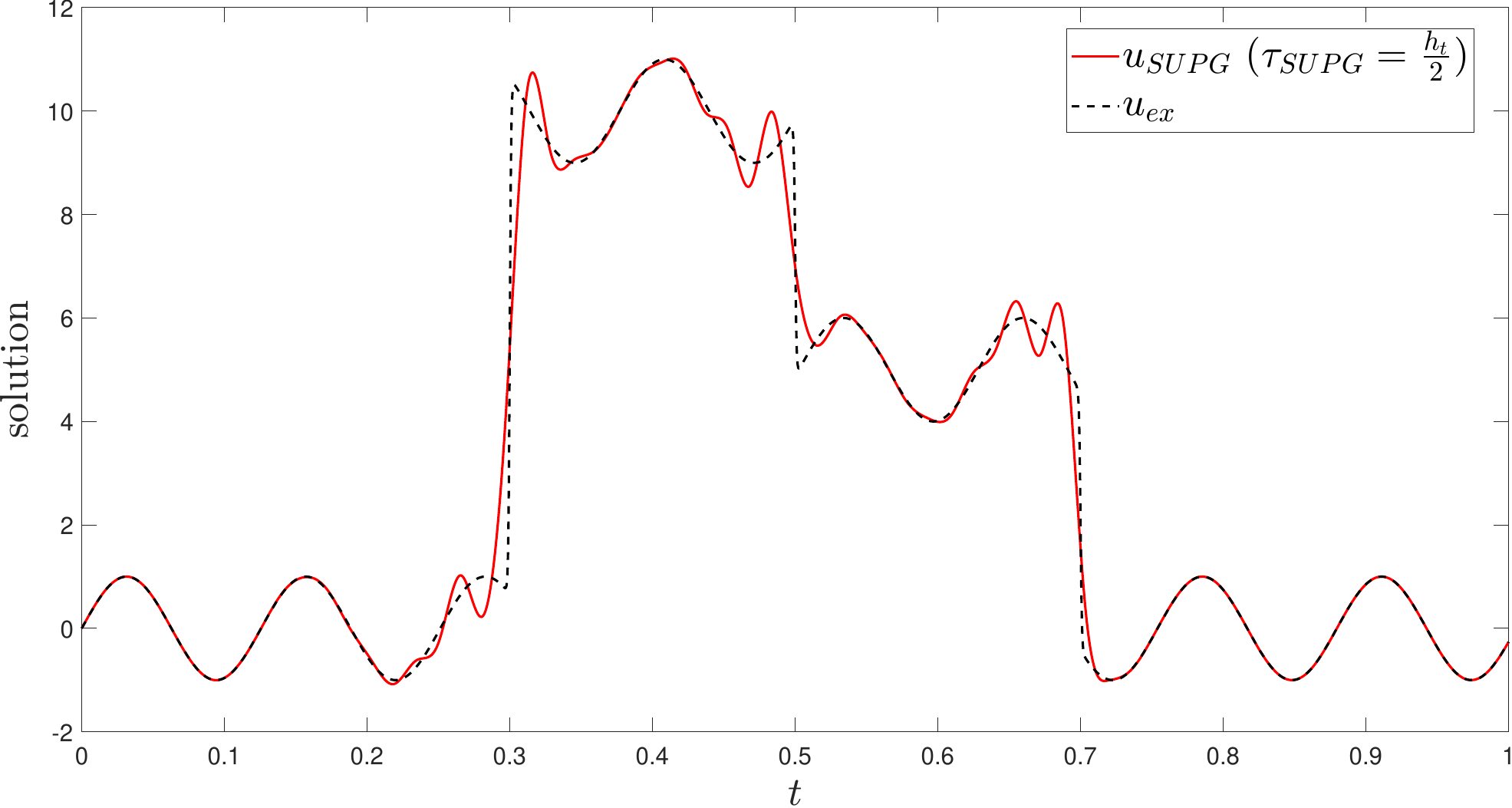}
	\caption{Advection equation, exact and SUPG solutions, with $h_{t}=2^{-6}$ and $p=3$.}
	\label{fig:SUPG_layers}
\end{figure}

\begin{figure}[htbp]
	\centering
		\includegraphics[width=1.00\textwidth]{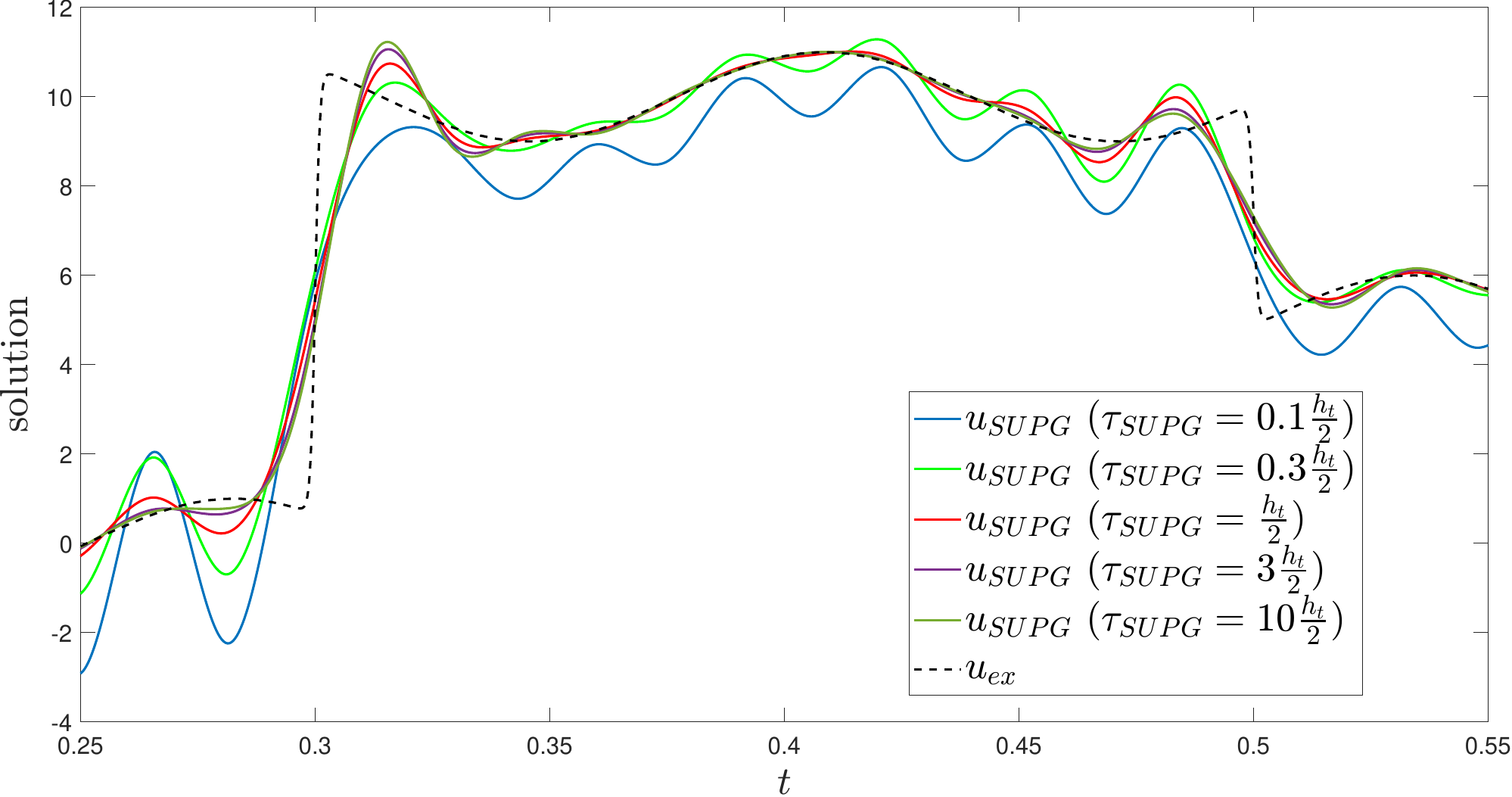}
	\caption{Advection equation, exact and SUPG solutions for different choices for $\tau_{\text{SUPG}}$, with $h_{t}=2^{-6}$ and $p=3$.}
	\label{fig:SUPG_tau}
\end{figure}

We select $f$ such that the exact solution is 
\begin{equation*} \label{sol_layers}
	u_{\text{ex}}(t)=\sin(50t)+10\frac{1+\tanh(\frac{t-t_0}{\delta})}{2}-5\frac{1+\tanh(\frac{t-t_1}{\delta})}{2}-5\frac{1+\tanh(\frac{t-t_2}{\delta})}{2},
\end{equation*}
with $t_0=0.3$, $t_1=0.5$, $t_2=0.7$ and $\delta=10^{-3}$.
Also on uniform meshes, due to the presence of internal layers, the standard Galerkin solution
is unstable (see Figure \ref{fig:GALERKIN_layers}). For high-degree
splines, plain SUPG has spurious oscillations, for any value of the stability parameter $\tau_{\text{SUPG}}$, see Figures \ref{fig:SUPG_layers} and \ref{fig:SUPG_tau} for the case $p=3$. In particular, as we can see in Figure \ref{fig:SUPG_tau}, $\tau_{\text{SUPG}}=\frac{h_t}{2}$ is the best value not only for $p=1$ (as presented in section \ref{sec:standard_stab}) but also for high-degree splines.

\begin{figure}[htbp]
	\centering
		\includegraphics[width=1.00\textwidth]{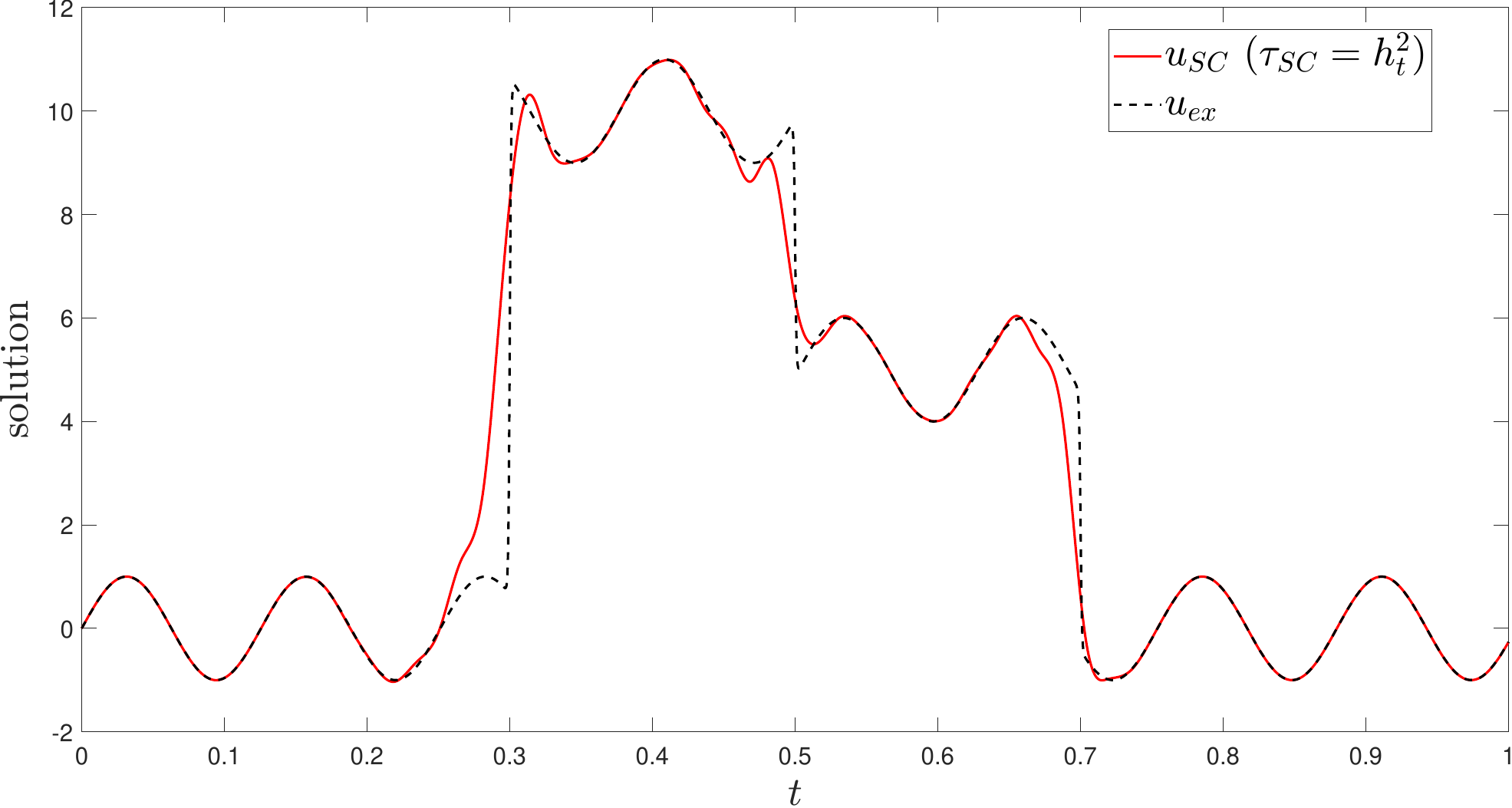}
	\caption{Advection equation, exact and Shock Capturing solutions, with $h_{t}=2^{-6}$ and $p=3$.}
	\label{fig:SC_layers}
\end{figure}

\begin{figure}[htbp]
	\centering
		\includegraphics[width=1.00\textwidth]{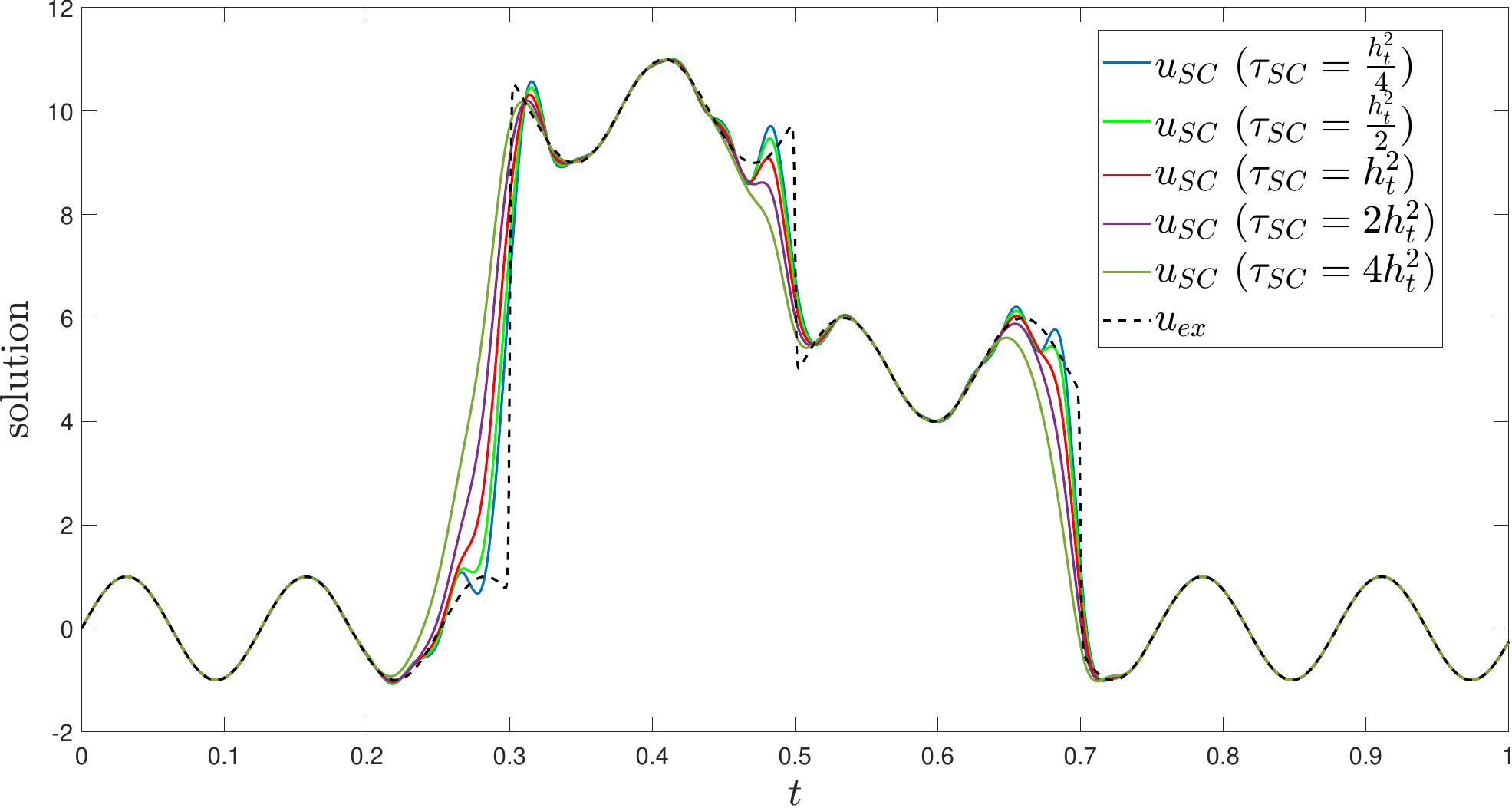}
 \caption{Advection equation, exact and Shock Capturing solutions for different choices for $\tau_{\text{SC}}$, with $h_{t}=2^{-6}$ and $p=3$.}
	\label{fig:SC_layers_tau}
\end{figure}

Adding Shock Capturing with $\tau_{\text{SC}}=h^2_t$ (Figure \ref{fig:SC_layers}) spurious oscillations are reduced but present for any $\tau_{\text{SC}}$ (Figure \ref{fig:SC_layers_tau}).

\begin{figure}[htbp]
	\centering
		\includegraphics[width=1.00\textwidth]{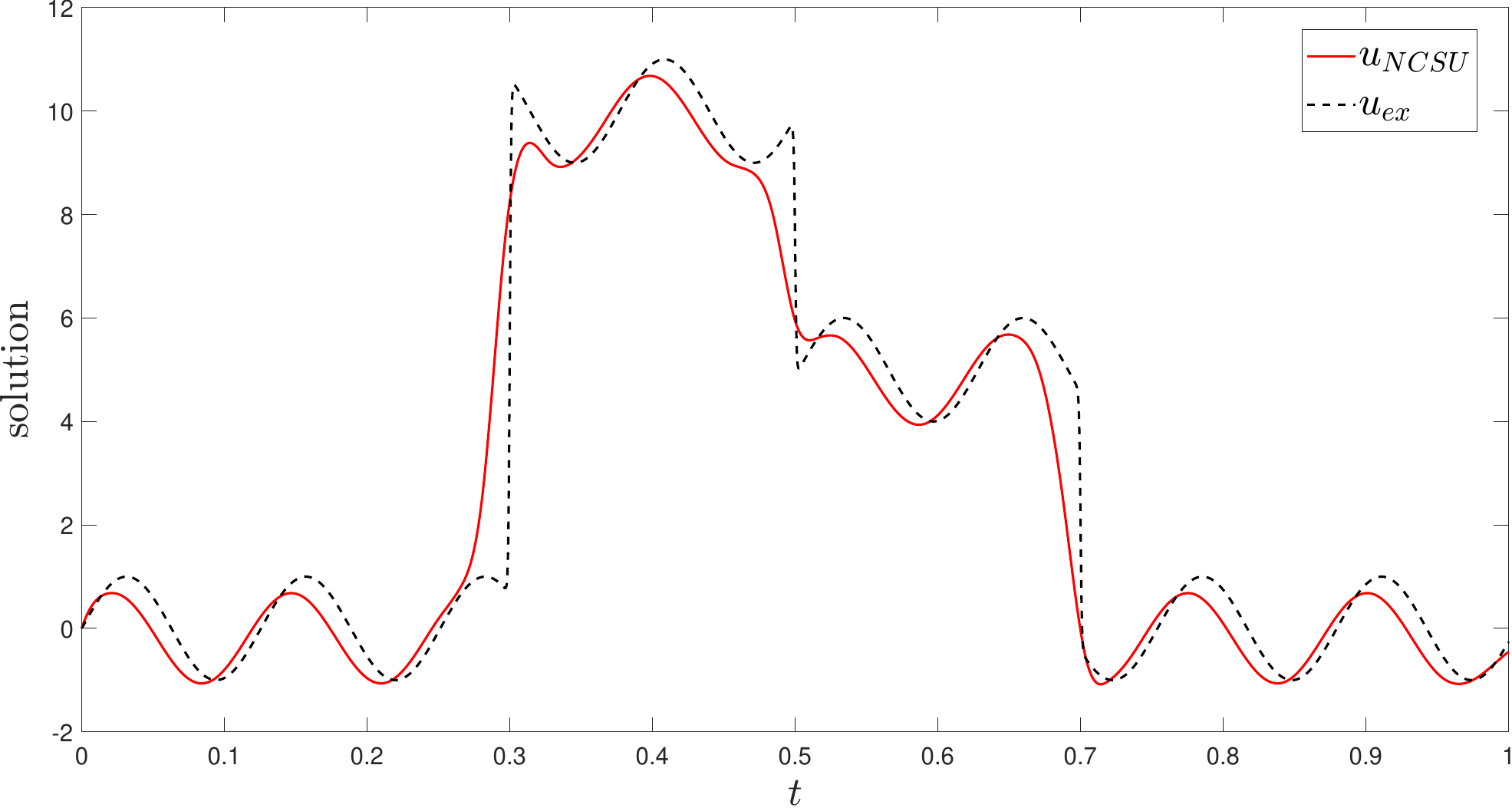}
	\caption{Advection equation, exact and NCSU solutions, with $h_{t}=2^{-6}$ and $p=3$.}
	\label{fig:NCSU_layers}
\end{figure}

With the non consistent NCSU method (Figure \ref{fig:NCSU_layers}) spurious oscillations disappear but the numerical and phase errors are significantly larger.

\begin{figure}[htbp]
	\centering
		\includegraphics[width=1.00\textwidth]{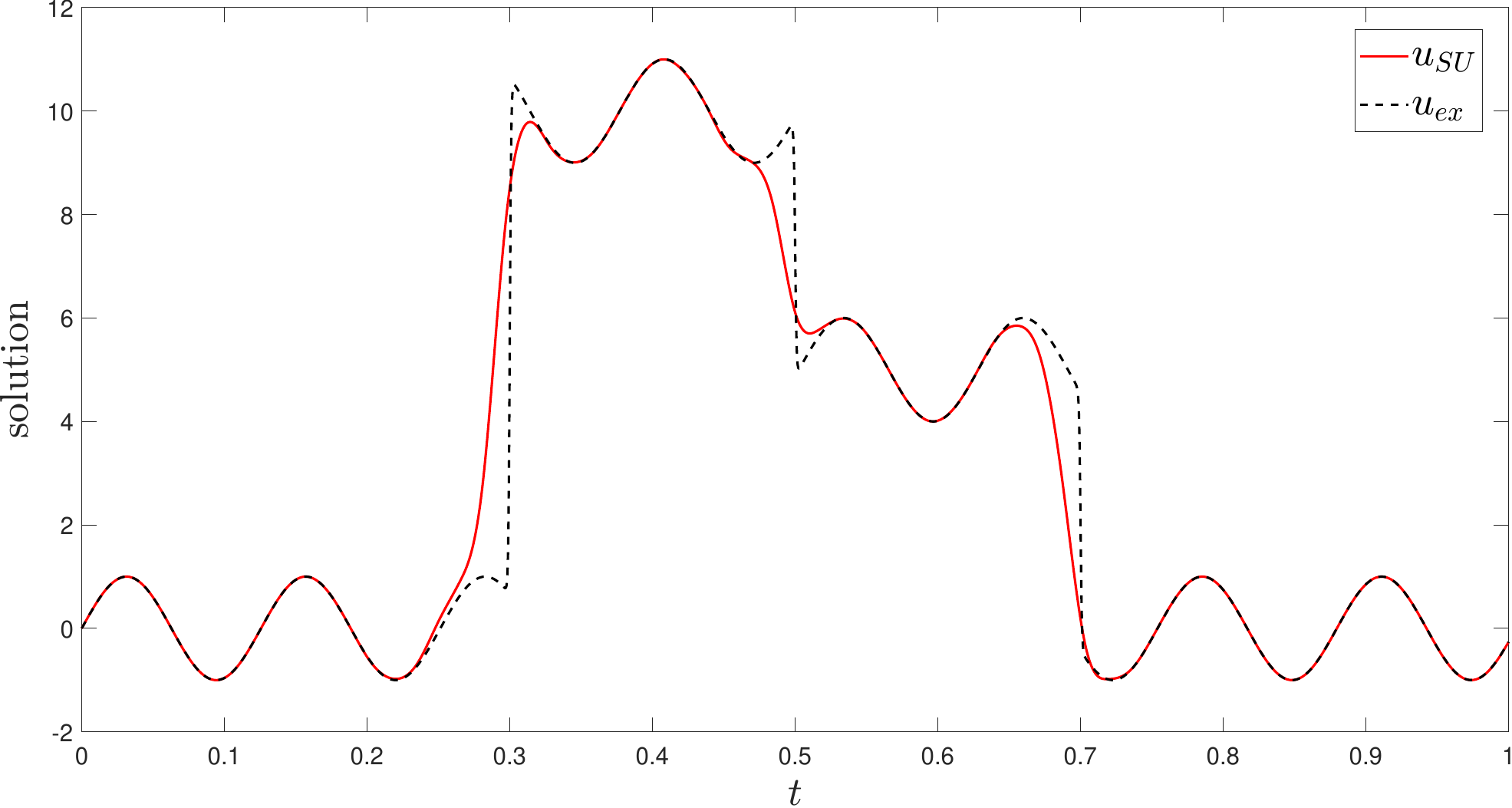}
	\caption{Advection equation, exact and SU solutions, with $h_{t}=2^{-6}$ and $p=3$.}
	\label{fig:SU_d3_layers}
\end{figure}

\begin{figure}[htbp]
	\centering
		\includegraphics[width=1.00\textwidth]{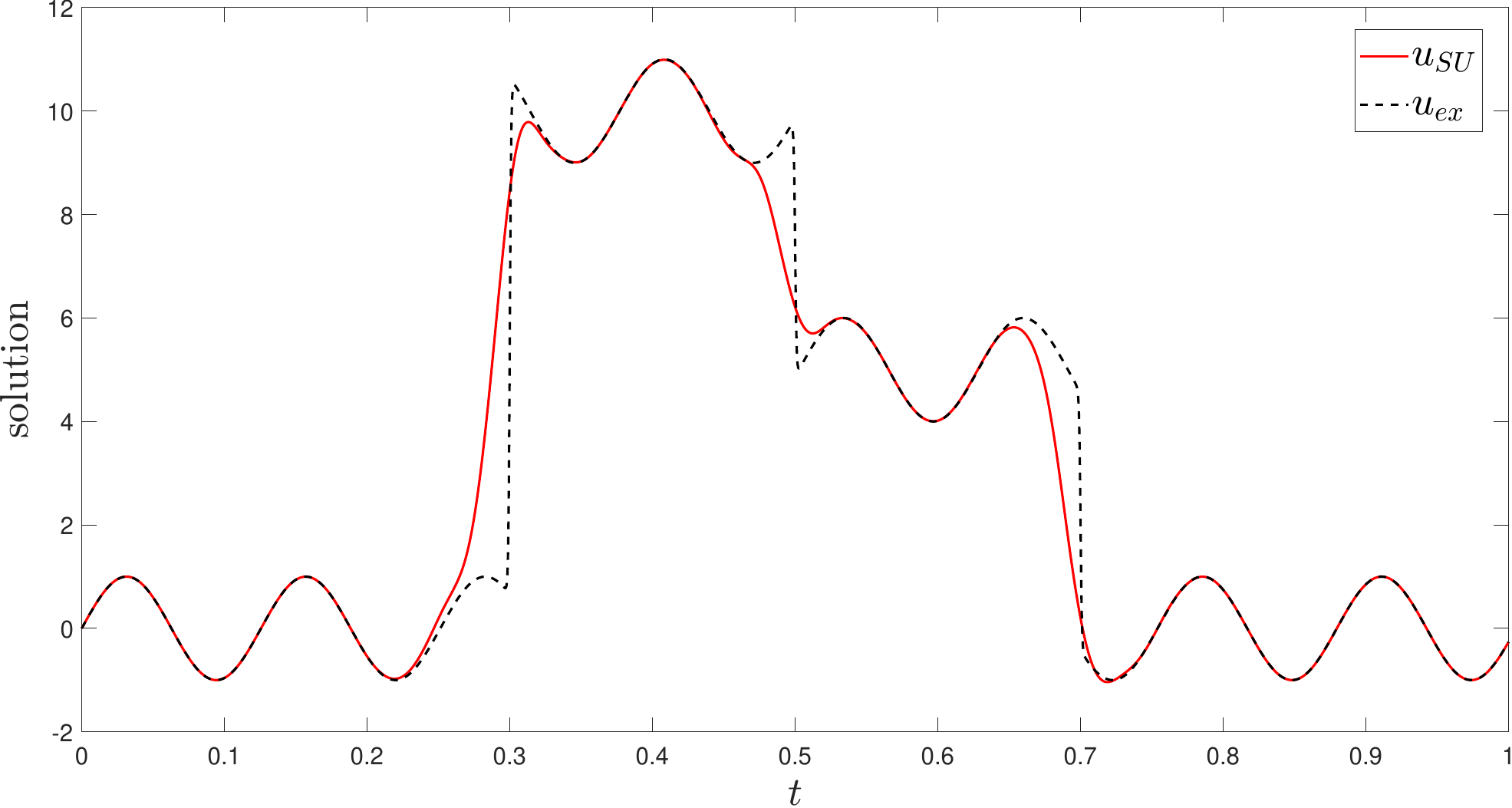}
	\caption{Advection equation, exact and SU solutions, with $h_{t}=2^{-6}$ and $p=4$.}
	\label{fig:SU_d4_layers}
\end{figure}

\begin{figure}[htbp]
	\centering
		\includegraphics[width=1.00\textwidth]{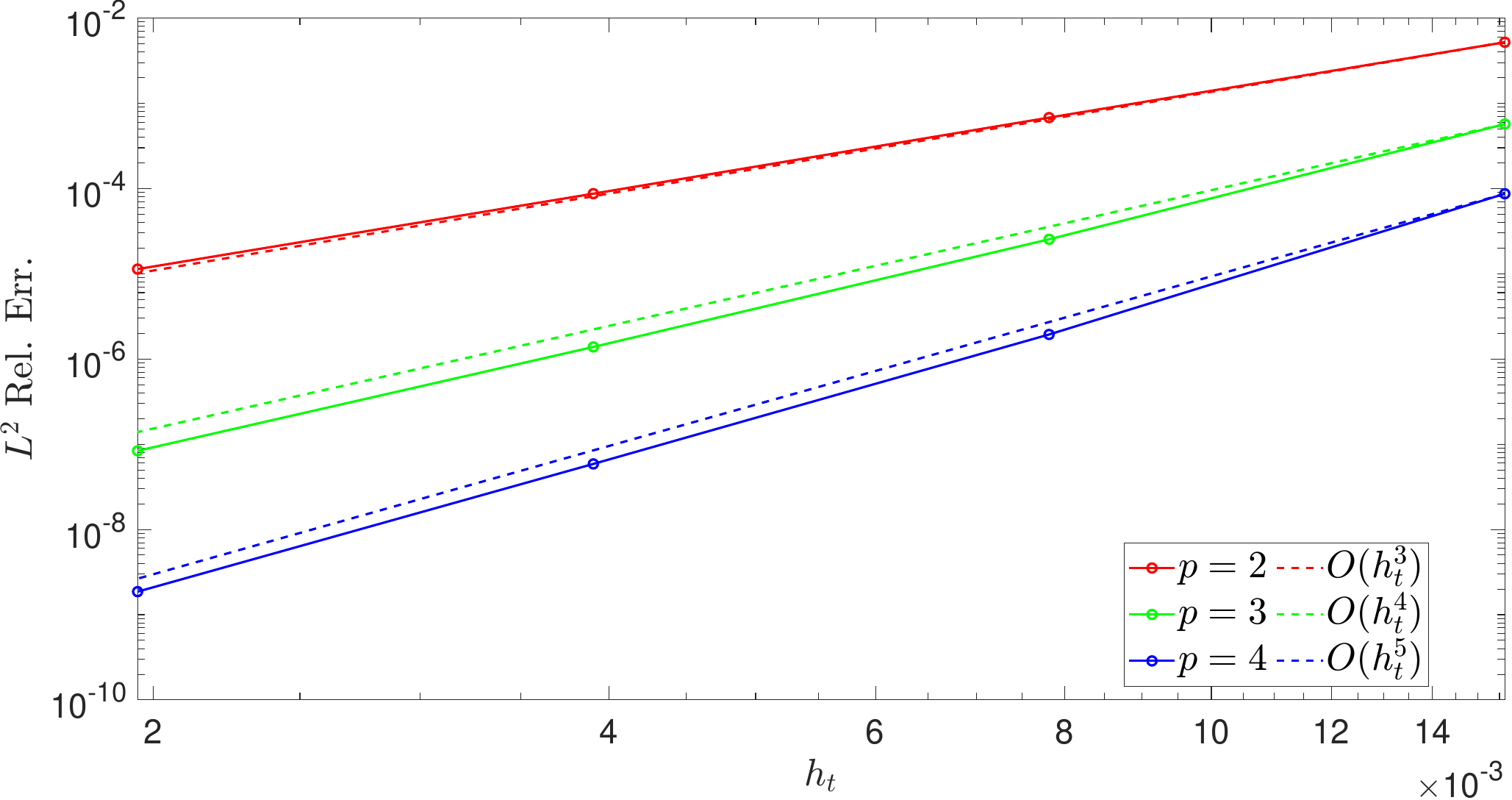}
	\caption{Advection equation, SU relative error plots in $L^2$-norm computed where the solution is smooth.}
	\label{fig:SU_L2_error_layers}
\end{figure}

Figures \ref{fig:SU_d3_layers} and \ref{fig:SU_d4_layers} show numerical results for $p=3$ and $p=4$, that assess the behavior of the SU formulation on uniform meshes: spurious oscillations are completely eliminated. Similar results are obtained from different degrees. 
Moreover relative error graphs in $L^2$-norm (Figure \ref{fig:SU_L2_error_layers}), calculated after the three layers where the solution is smooth ($t>0.85$), show that the relative error converges optimally. 

\begin{figure}[htbp]
	\centering
		\includegraphics[width=1.00\textwidth]{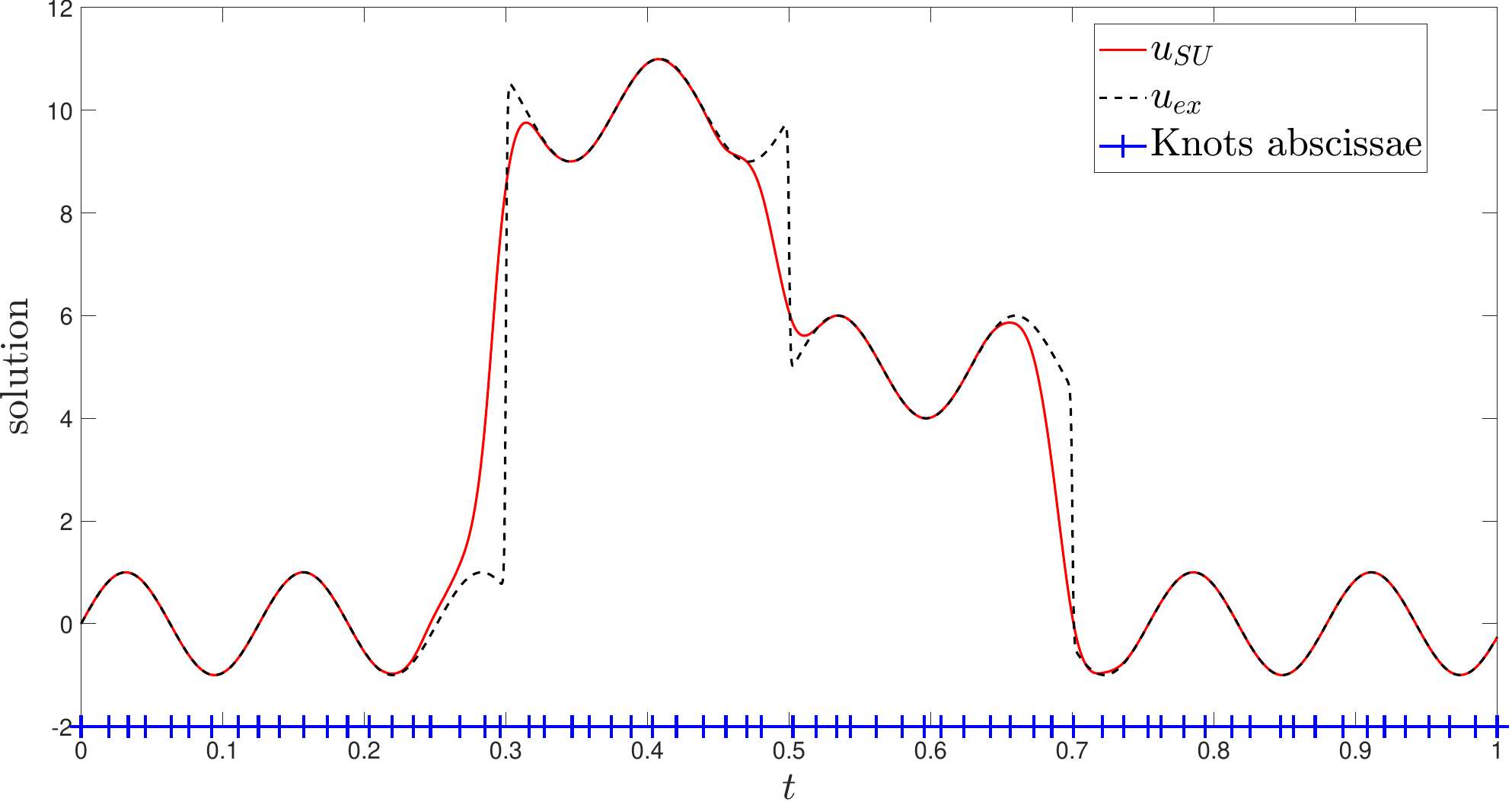}
	\caption{Advection equation, exact and SU solutions on a non-uniform mesh (depicted in blue on the horizontal axis) and $p=3$.}
	\label{fig:SU_d3_layers_nonunif}
\end{figure}

Stable and accurate behavior of SU method is also confirmed if we deal with non-uniform meshes, as we can see in Figure \ref{fig:SU_d3_layers_nonunif}.

\subsection{Advection-diffusion equation}

\begin{figure}[htbp]
	\centering
		\includegraphics[width=1.00\textwidth]{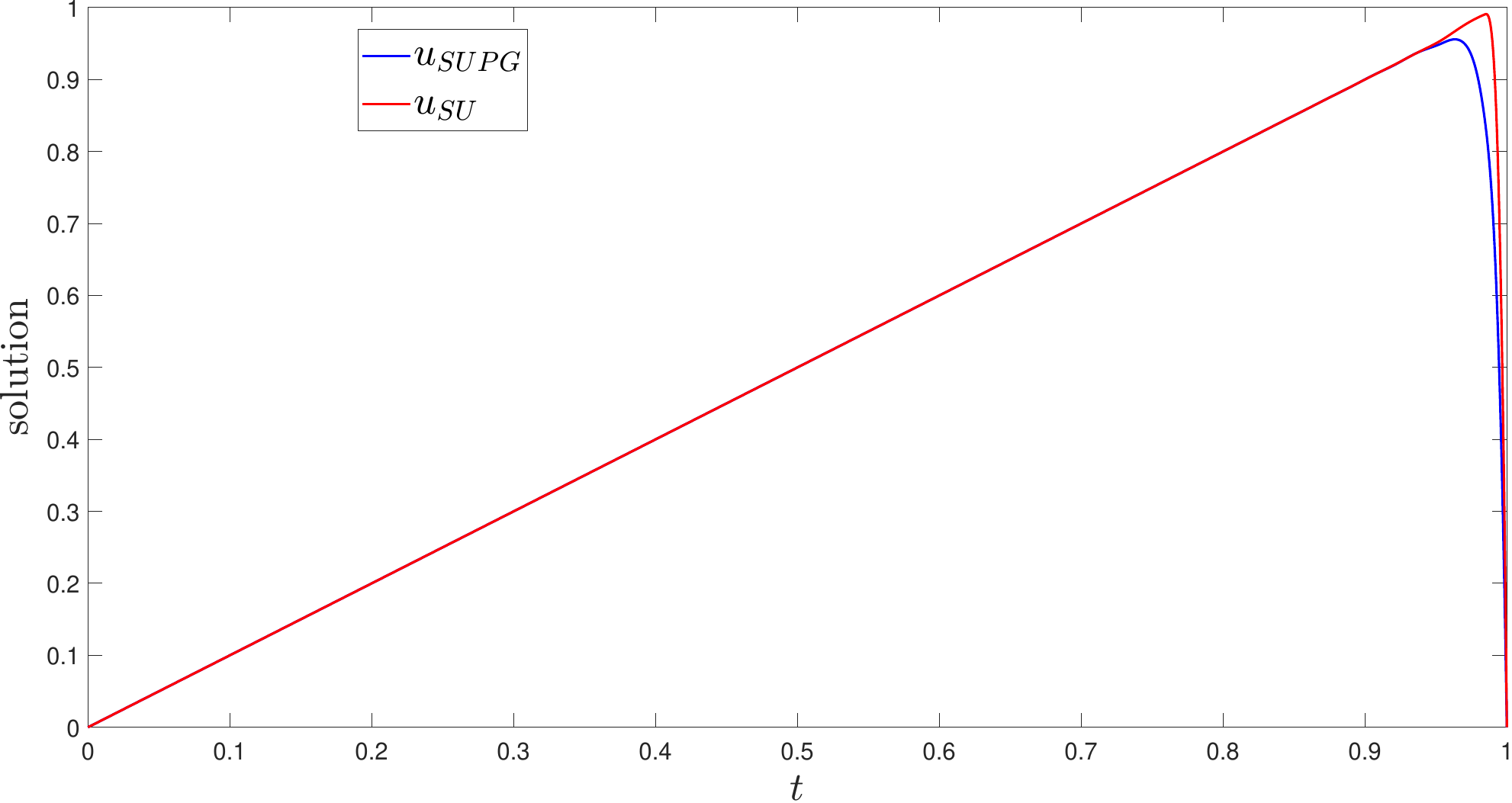}
	\caption{Advection-diffusion equation, SUPG and SU solution; $h_{t}=2^{-6}$ and $p=3$.}
	\label{fig:adv_diff_d3}
\end{figure}

\begin{figure}[htbp]
	\centering
		\includegraphics[width=1.00\textwidth]{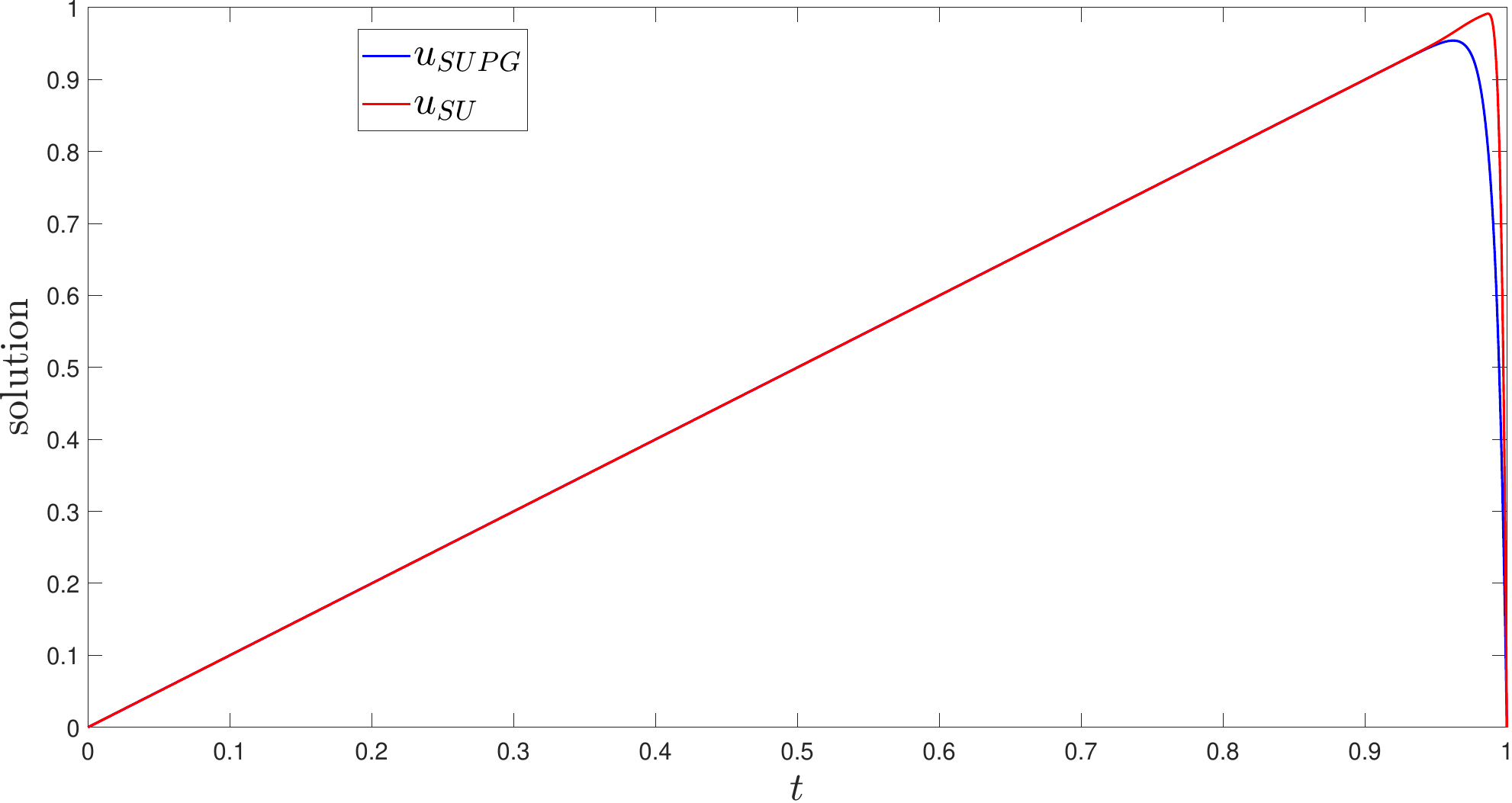}
	\caption{Advection-diffusion equation, SUPG and SU solution; $h_{t}=2^{-6}$ and $p=4$.}
	\label{fig:adv_diff_d4}
\end{figure}

As in Section \ref{sec:adv_results} we consider a uniform mesh.
We consider the advection-diffusion equation \eqref{eq:adv_diff_problem} on $(0,T)$ with $T=1$, $f = 1$ and $\varepsilon=10^{-6}$.
Figures \ref{fig:adv_diff_d3} and \ref{fig:adv_diff_d4} show SUPG and
SU solutions.The SU method demonstrates higher
  accuracy, albeit at a higher computational cost. In our
  implementation the  increased cost arises from the need to evaluate
  the residual at each fixed point iteration. However, the issue of computational cost and efficient implementation deserves a more in-depth exploration beyond the scope of this work.

\subsection{Heat equation} \label{sec:heat_num_setting}

\subsubsection{2D space-time domain}\label{sec:heat_2D}

For the first test we consider $\Omega=(0,1)$, we solve a uniform mesh and we use the same mesh-size in space and in time, i.e. we  set $h_s=h_t=:h$.
We consider the heat equation \eqref{eq:heat_problem} on $(0,1)\times(0,T)$ with $T=1$ and $f$ as follows:
\begin{equation*} 
	f(x,t)= \delta^{-2} \exp(-((x - (1/4(\sin(10\pi t) + 2)))/\delta)^2) \ \chi_{[0.3,0.6]}(t) ,
\end{equation*}
where
\begin{equation*} 
\chi_{[0.3,0.6]}(t)= 
	\left\{
	\begin{array}{llllll}
		  1 & \ \ \text{for} \ t \in [0.3,0.6],\\
		  0 &  \ \ \text{otherwise},
	\end{array}
	\right.
\end{equation*}
and with $\delta=10^{-3}$.

\begin{figure}[htbp]
	\centering
			\includegraphics[width=1\textwidth]{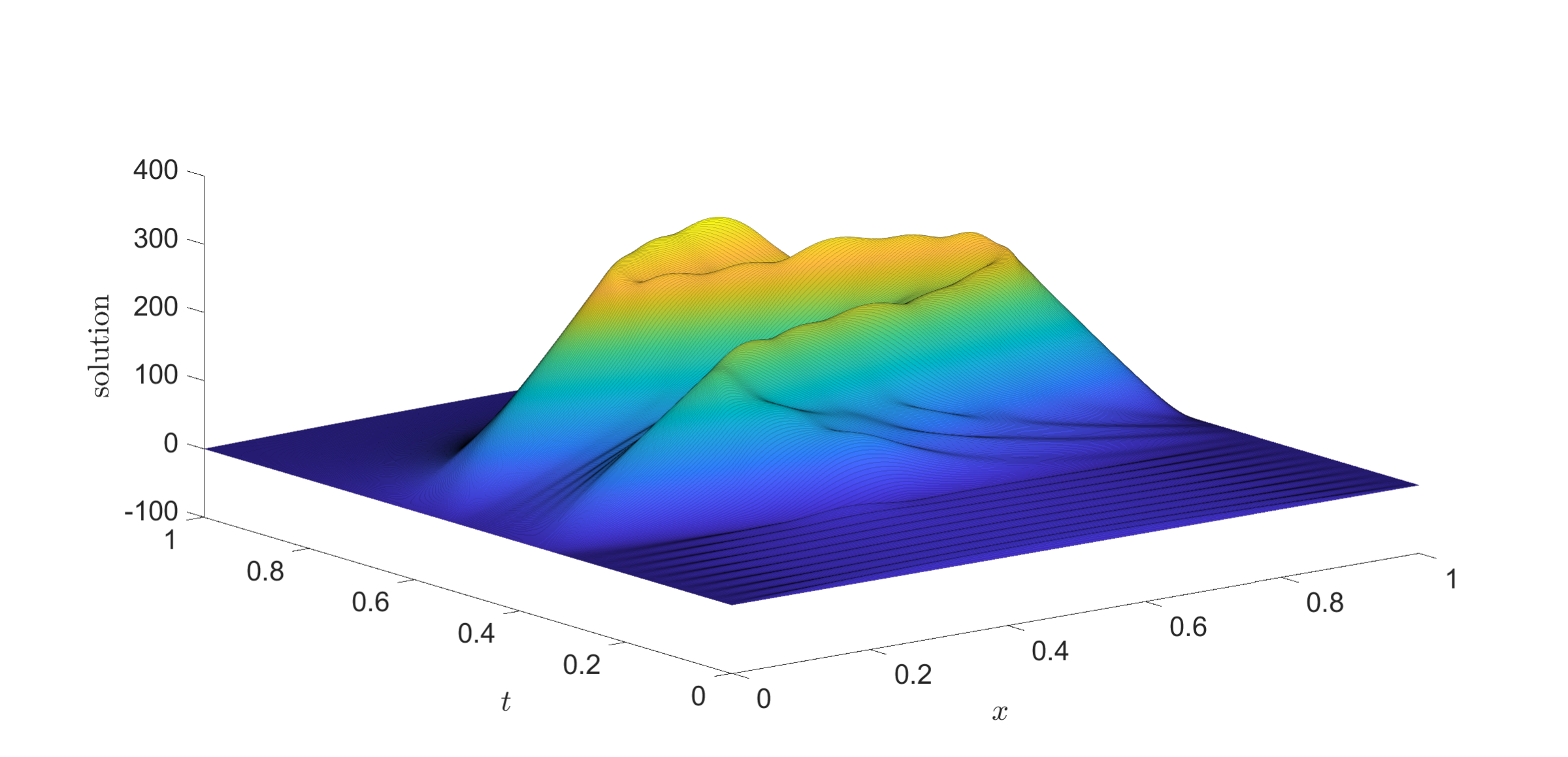}
	\caption{Numerical solution of the Galerkin method for the heat equation ($\Omega=(0,1)$), with $h=2^{-6}$ and $p=3$, spurious oscillations are present also for $t<0.3$.}
	\label{fig:ST_segment_GALERKIN}
\end{figure}

\begin{figure}[htbp]
	\centering
			\includegraphics[width=1\textwidth]{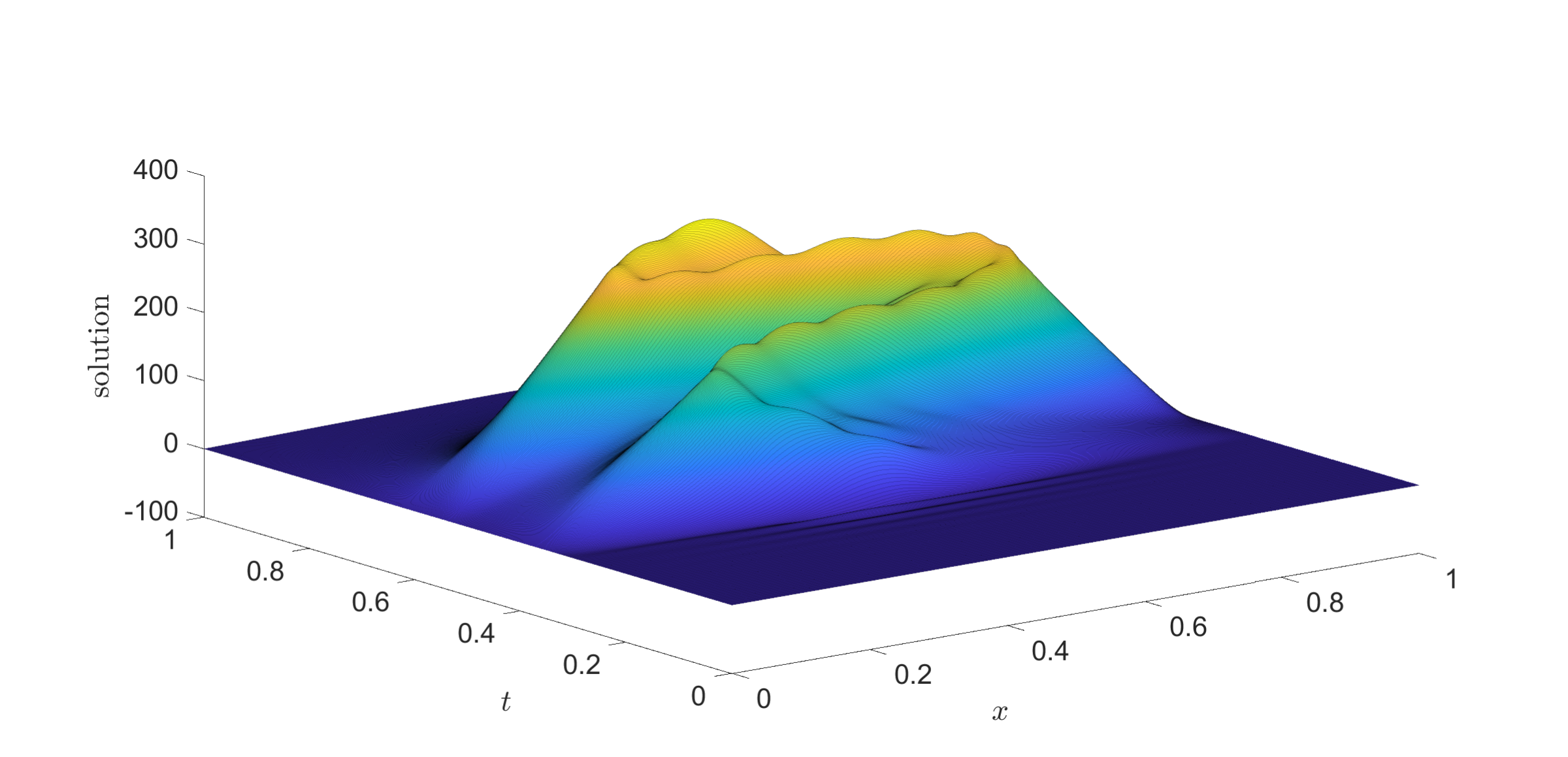}
	\caption{Numerical solution of SUPG method for the heat equation ($\Omega=(0,1)$), with $h=2^{-6}$ and $p=3$.}
	\label{fig:ST_segment_SUPG}
\end{figure}

\begin{figure}[htbp]
	\centering
			\includegraphics[width=\textwidth]{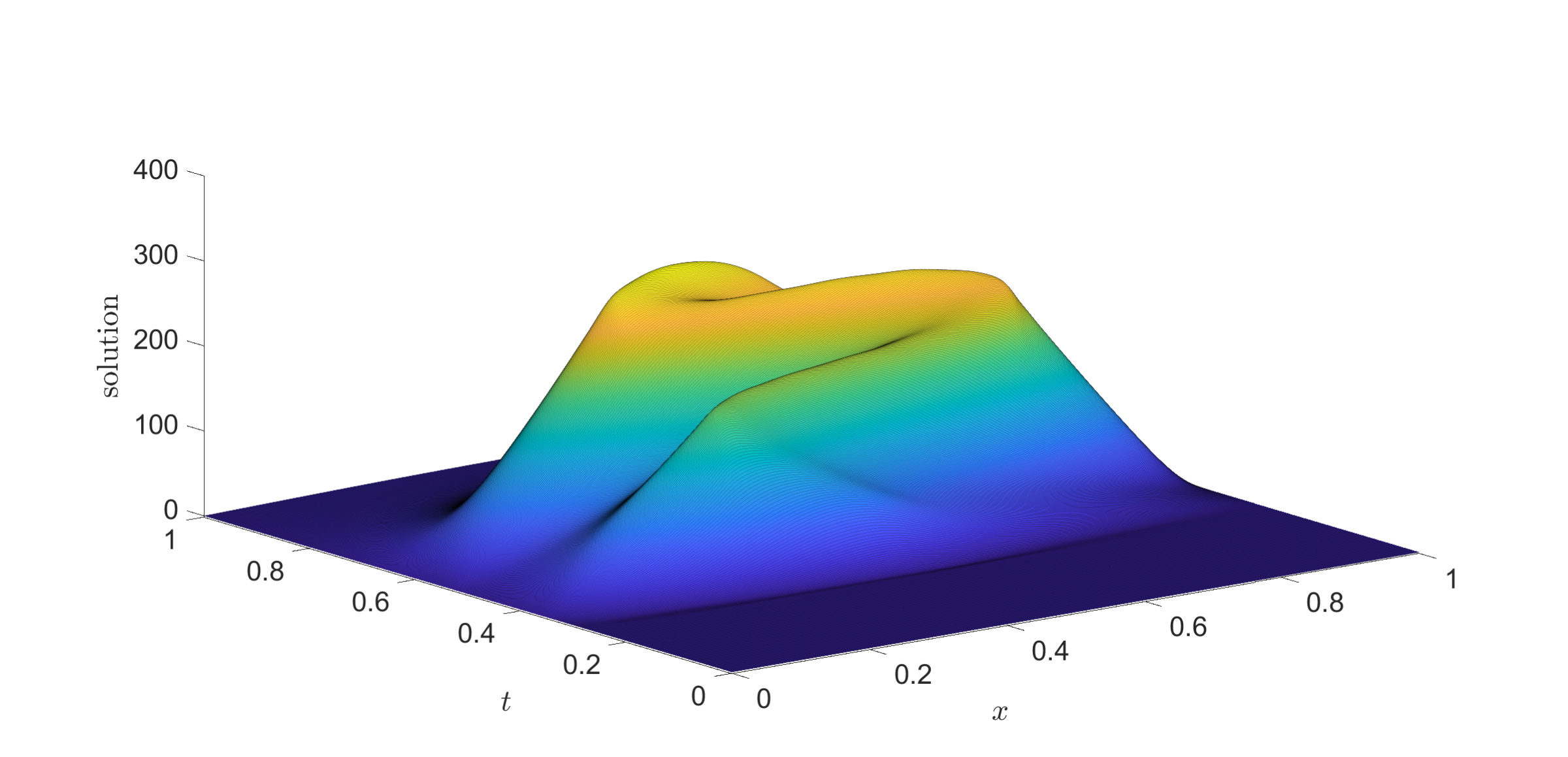}
	\caption{Numerical solution of the SU method for the heat equation ($\Omega=(0,1)$), with $h=2^{-6}$ and $p=3$.}
	\label{fig:ST_segment_SU}
\end{figure}

\begin{figure}[htbp]
	\centering
			\includegraphics[width=1\textwidth]{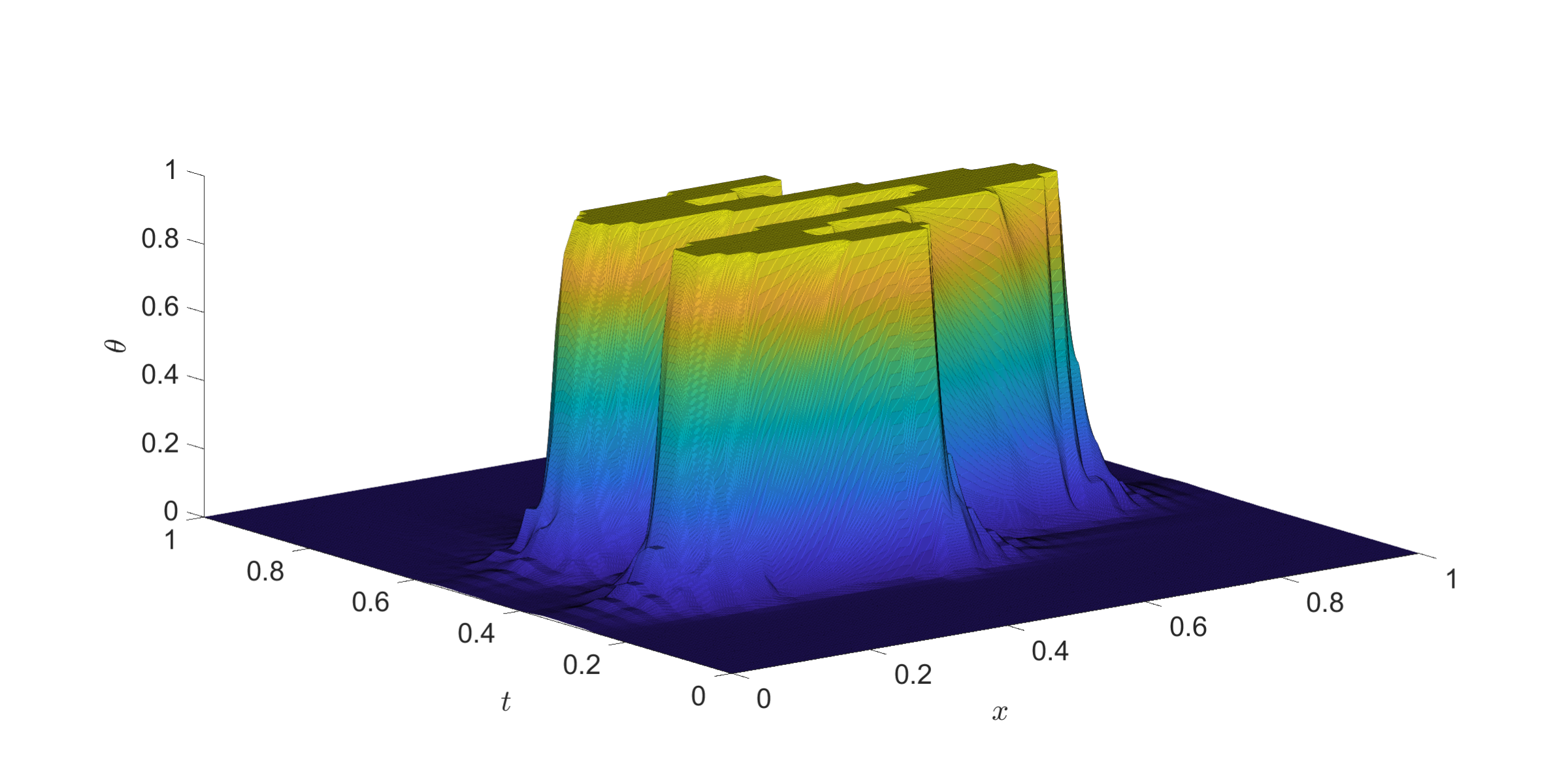}
	\caption{function $\theta(x,t)$ ($\Omega=(0,1)$), with $h=2^{-6}$ and $p=3$.}
	\label{fig:ST_segment_SU_theta}
\end{figure}

In Figures \ref{fig:ST_segment_GALERKIN}, \ref{fig:ST_segment_SUPG}
and  \ref{fig:ST_segment_SU} the numerical solutions by
Galerkin, SUPG and the SU methods are presented. 
Figure \ref{fig:ST_segment_SU_theta} displays the graph of the function $\theta(x,t)$, which plays a crucial role in activating the high-order Upwind stabilization in the proximity of the layers. The function $\theta(x,t)$ serves as a key indicator, guiding the activation of the stabilization technique to effectively address the presence of sharp layers in the solution.

\subsubsection{3D space-time domain}

\begin{figure}[htbp]
	\centering
		\includegraphics[width=0.4\textwidth]{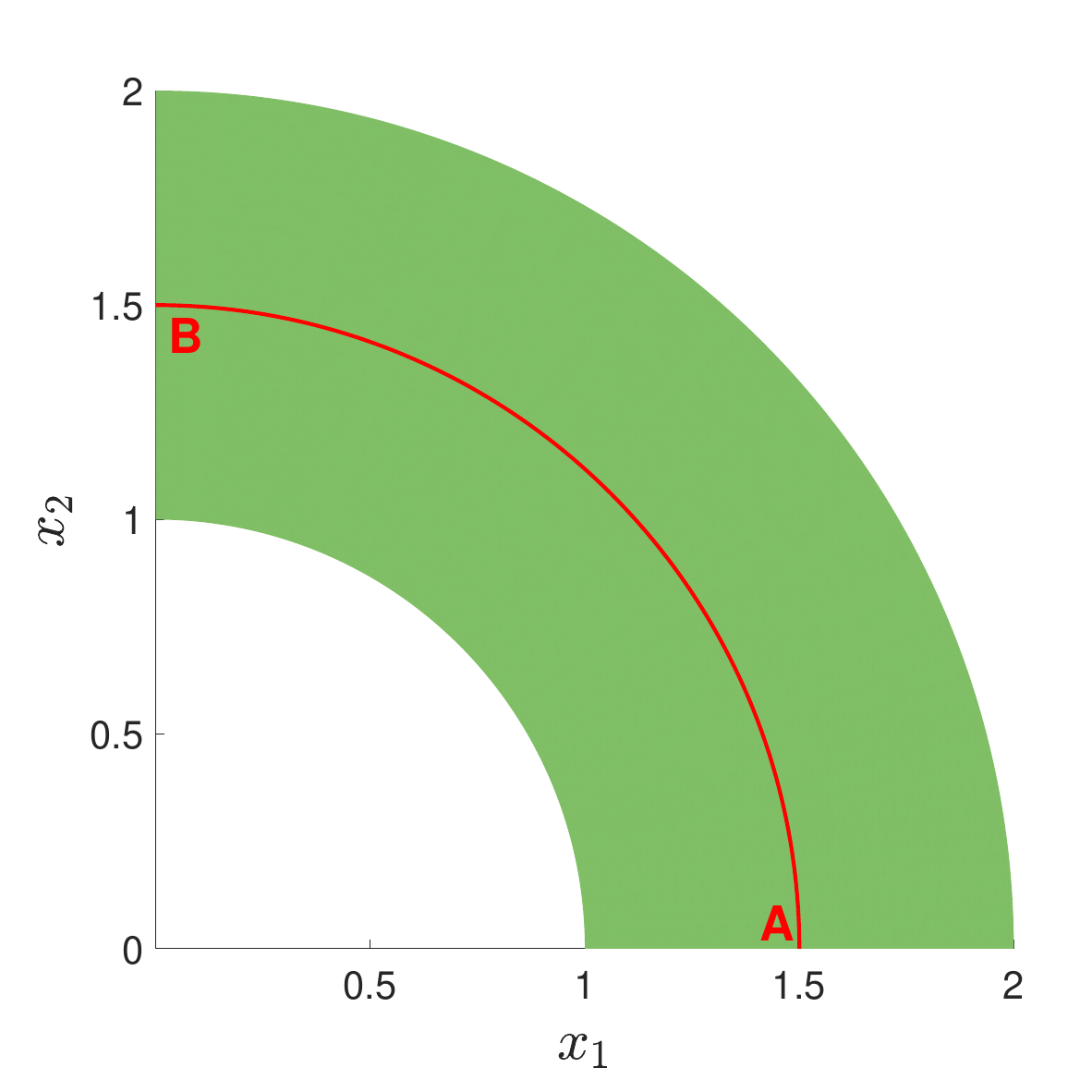}
	\caption{Quarter annulus with section line A-B.}
	\label{fig:qannulus}
\end{figure}

For the second test we take as $\Omega$ a quarter annulus (Figure \ref{fig:qannulus}), and set
\begin{equation*} 
	f(x_1,x_2,t)= 10^3/(2\pi\delta^2)\exp(-1/2(((x_1-1.5\cos(\pi/2t))/\delta)^2+((x_2-1.5\sin(\pi/2t))/\delta)^2))\ \chi_{[0.3,0.6]}(t),
\end{equation*}
with $\delta=0.1$.

\begin{figure}[htbp]
	\centering
			\includegraphics[width=1\textwidth]{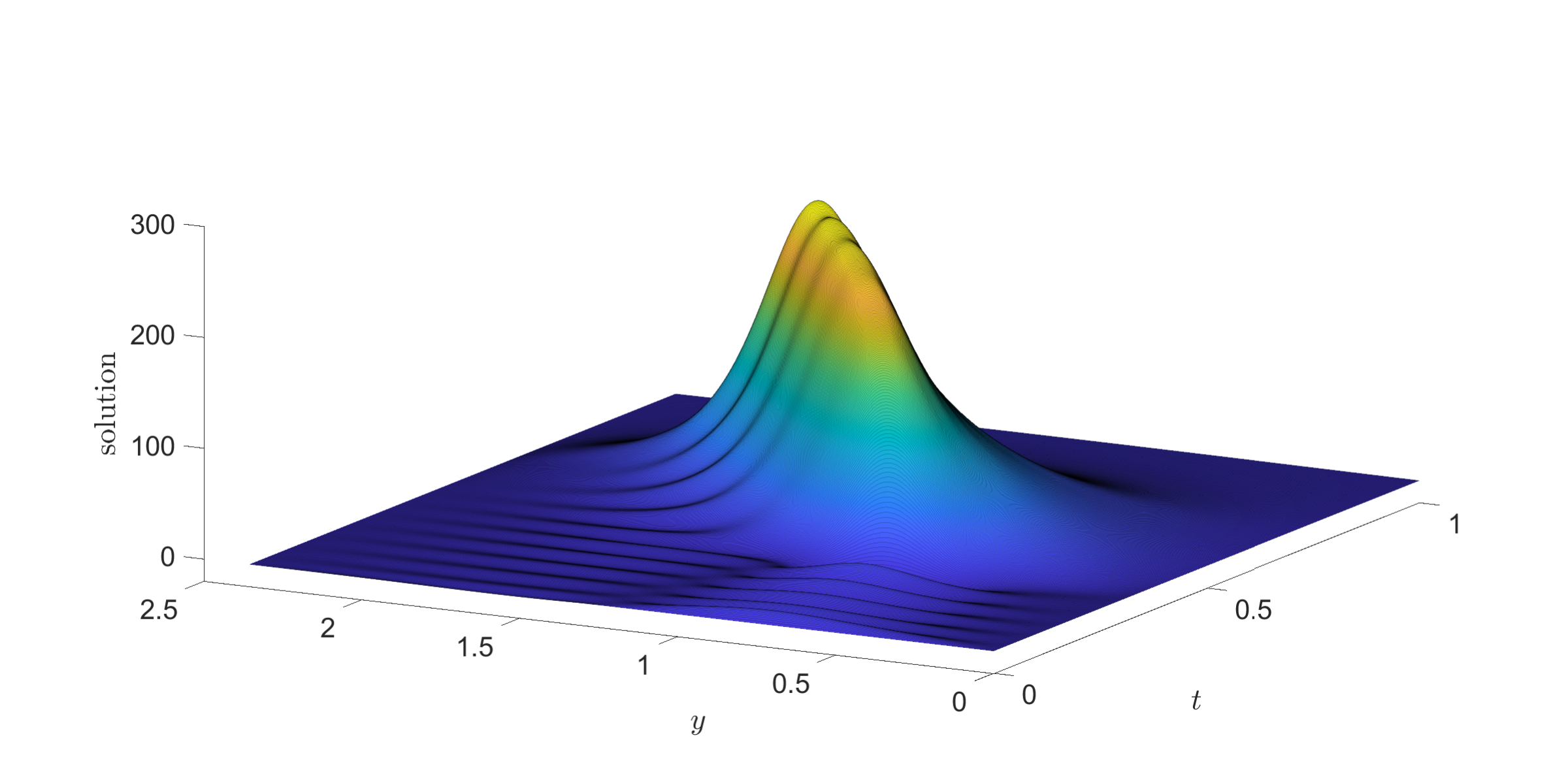}
	\caption{Numerical solution of Galerkin method  (along section A-B) for the heat equation on the quarter annulus, with $h=2^{-5}$ and $p=3$, spurious oscillations are present also for $t<0.3$.}
	\label{fig:ST_qannulus_GALERKIN}
\end{figure}

\begin{figure}[htbp]
	\centering
			\includegraphics[width=1\textwidth]{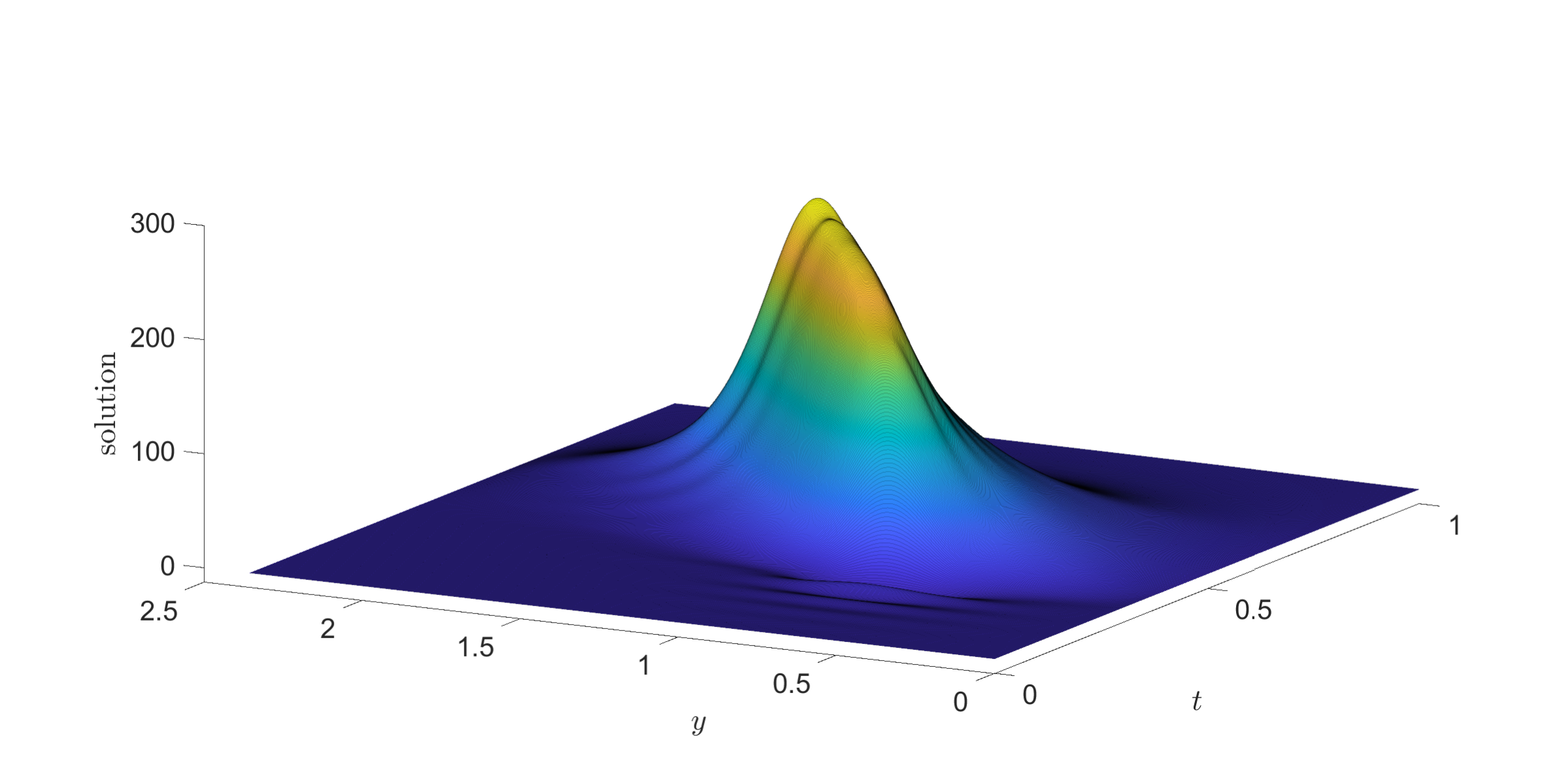}
	\caption{Numerical solution of SUPG method (along section A-B) for the heat equation on the quarter annulus, with $h=2^{-5}$ and $p=3$.}
	\label{fig:ST_qannulus_SUPG}
\end{figure}

\begin{figure}[htbp]
	\centering
			\includegraphics[width=1\textwidth]{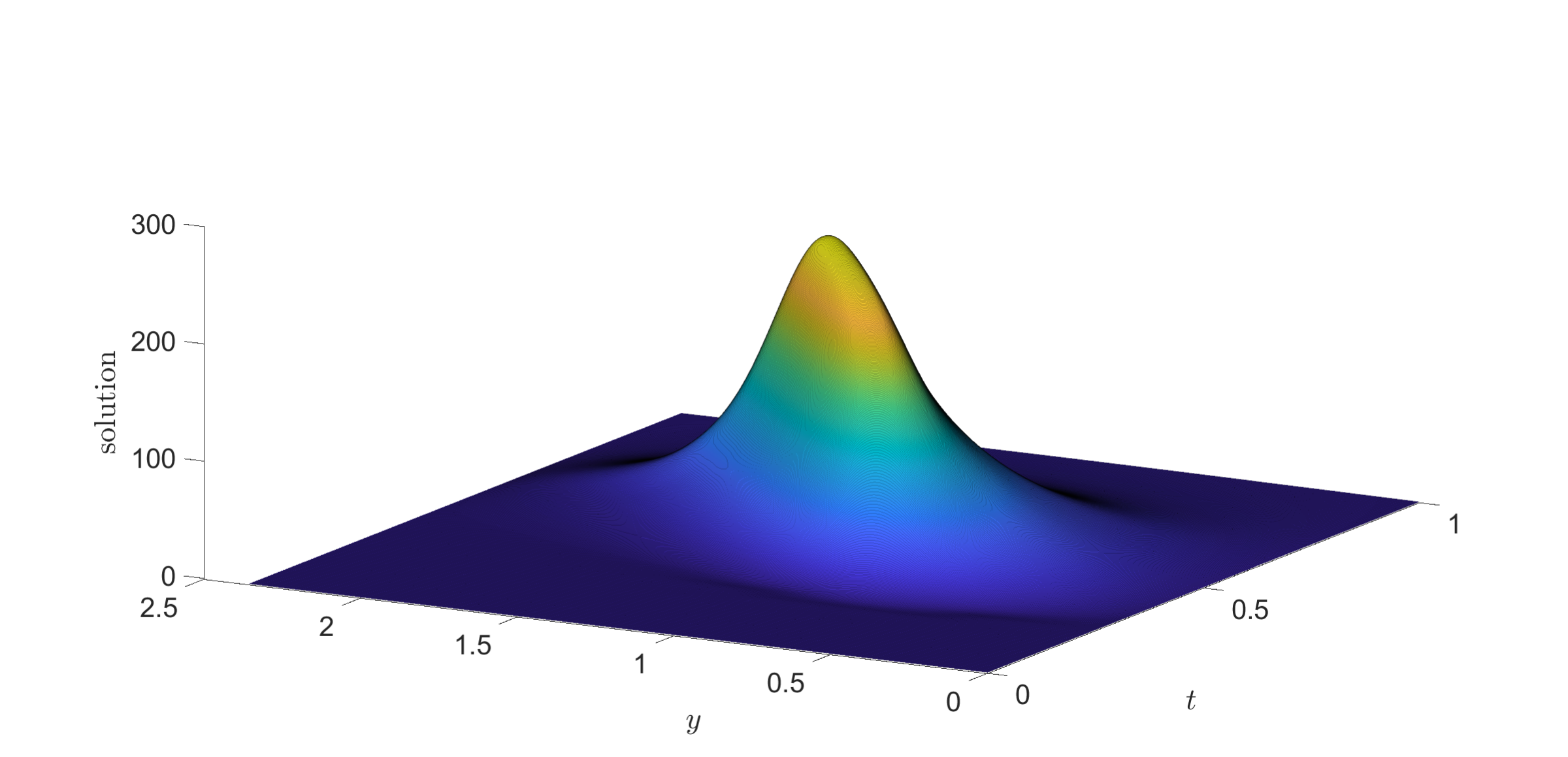}
	\caption{Numerical solution of the SU method (along section A-B) for the heat equation on the quarter annulus, with $h=2^{-5}$ and $p=3$.}
	\label{fig:ST_qannulus_SU}
\end{figure}

In Figures \ref{fig:ST_qannulus_GALERKIN} and
\ref{fig:ST_qannulus_SUPG} the numerical results that assess the
behavior of the space-time Galerkin approximation and SUPG method are
presented. As in Section \ref{sec:heat_2D}, we observe the emergence of spurious oscillations, particularly in the case of the plain Galerkin method. However, when examining the numerical results of the SU method in Figure \ref{fig:ST_qannulus_SU}, we can see that they are free from spurious oscillations.

\section{Conclusions} 

In this work, we have presented a novel
  space-time method for the heat equation in the framework of IgA.  It
  is based on smooth spline approximation in time and incorporates a
  stabilizing term that extends the SUPG stabilization mechanism to
  high-degree and continuity splines, promoting  causality with
  respect to time.

We have conducted various numerical benchmarks to validate our method.
 The results provided numerical
  evidence of the optimal order of convergence for smooth solutions
  and of the stable behavior even in the presence of sharp layers and concentrated source terms.

While our focus in this work was not on
  computational cost, we acknowledge the significance of efficient and
  fast solvers in the space-time framework. In particular the higher dimensionality poses computational challenges that need to be addressed. We plan to explore and address these computational aspects in future works.


\section*{Acknowledgements}
The authors are members of the Gruppo Nazionale Calcolo Scientifico-Istituto
Nazionale di Alta Matematica (GNCS-INDAM). 
Support for this research was partially provided by a grant through Regione Lombardia, POR FESR 2014-2020 - Call HUB Ricerca e Innovazione, Progetto 1139857 CE4WE: Approvvigionamento energetico e gestione della risorsa idrica nell’ottica dell’Economia Circolare (Circular Economy for Water and Energy). G. Loli was also partially supported by the GNCS-INdAM through the ``Bando Finanziamento Giovani Ricercatori 2021-2022 GNCS''.
The authors acknowledge the contribution of the National Recovery and Resilience Plan, Mission 4 Component 2 - Investment 1.4 - NATIONAL CENTER FOR HPC, BIG DATA AND QUANTUM COMPUTING, spoke 6.

\bibliographystyle{elsarticle-num}
\bibliography{SU_method}

\end{document}